\documentclass[12pt]{article}
\usepackage[top=25.4mm, bottom=25.4mm, left=25.4mm, right=25.4mm]{geometry}

\usepackage{amssymb,amstext,amsmath}
\usepackage{graphicx}% Include figure files
\usepackage{dcolumn}% Align table columns on decimal point
\usepackage{bm} % bold math
\usepackage{multirow}
\usepackage{bigstrut}
\usepackage{makecell,rotating}
\usepackage{mathrsfs}
\usepackage{booktabs}
\usepackage{threeparttable}
\usepackage{chngpage}
\usepackage{color}
\usepackage{caption,subcaption}
\usepackage{cite}
\usepackage{xr}
\newcommand{\1}{\mathbf{x}}
\newcommand{\A}{\mathbf{A}}
\newcommand{\B}{\mathbf{B}}
\newcommand{\C}{\mathbf{C}}
\newcommand{\W}{\mathbf{W}}
\newcommand{\U}{\mathbf{u}}
\newcommand{\y}{\mathbf{y}}
\newcommand{\R}{\mathbb{R}}
\newcommand{\N}{\mathbb{N}}
\newcommand{\T}{\text{T}}
\newcommand{\2}{\mathbf{R}}
\newcommand{\M}{\mathbf{M}}
\newcommand{\h}{\mathbf{H}}

\newcommand{\p}{\mathbf{P}}
\newcommand{\la}{\mathbf{\Lambda}}
\newcommand{\Q}{\mathbf{Q}}
\newcommand{\I}{\mathbf{I}}
\newcommand{\e}{\mathbf{e}}
\newcommand{\X}{\mathbf{\Theta}}
\newcommand{\0}{\mathbf{0}}

\usepackage{color}
	 \definecolor{darkred}{rgb}{0.75,0,0}
	 \definecolor{darkgreen}{rgb}{0,0.5,0}
	 \definecolor{darkblue}{rgb}{0,0,0.75}
  	 \definecolor{darkorange}{rgb}{1,0.9,0.1}
	 \definecolor{dark}{rgb}{0,0,0}

\begin{document}

\title{Energy cost for target control of complex networks}
%\title{Minimum energy control in complex networks}
\author{Gaopeng Duan$^{1^*}$, Aming Li$^{2,3^*}$, Tao Meng$^{1}$, and Long Wang$^{1^\dag}$}
\date{\today}
\maketitle
\begin{enumerate}
  \item Center for Systems and Control, College of Engineering, Peking University, Beijing 100871, China
   \item Department of Zoology, University of Oxford, Oxford OX1 3PS, UK
      \item Department of Biochemistry, University of Oxford, Oxford OX1 3QU, UK
   \item[$*$]  These authors contributed equally to this work
   \item[$\dag$] Correspondence to: longwang@pku.edu.cn
\end{enumerate}

\begin{abstract}
To promote the implementation of realistic control over various complex networks, recent work has been focusing on analyzing energy cost.
Indeed, the energy cost quantifies how much effort is required to drive the system from one state to another when it is fully controllable.
A fully controllable system means that the system can be driven by external inputs from any initial state to any final state in finite time.
However, it is prohibitively expensive and unnecessary to confine that the system is fully controllable when we merely need to accomplish the so-called target control---controlling a subnet of nodes chosen from the entire network.
Yet, when the system is partially controllable, the associated energy cost remains elusive.
Here we present the minimum energy cost for controlling an arbitrary subset of nodes of a network.
Moreover, we systematically show the scaling behavior of the precise upper and lower bounds of the minimum energy in term of the time given to accomplish control.
For controlling a given number of target nodes, we further demonstrate  that the associated energy over different configurations can differ by several orders of magnitude.
When the adjacency matrix of the network is nonsingular, we can simplify the framework by just considering the induced subgraph spanned by target nodes instead of the entire network.
Importantly, we find that, energy cost could  be saved by orders of magnitude as we only need the partial controllability of the entire network.
Our theoretical results are all corroborated by numerical calculations, and pave the way for estimating the energy cost to implement realistic target control in various applications.
\end{abstract}

\section{Introduction}

Network control has received much attention in the past decade
\cite{barabasi2016network,Havlin2004book,duan2017asynchronous,Liu2016Rev,liu2013observability,wang2007new,wang2007finite}.
The practical requirement of controlling complex networks from arbitrary initial to final states in finite time by appropriate external inputs motivates various explorations on the essential attribute of complex networked systems---network controllability \cite{chen2017pinning,Duan2019PRE,guan2017controllability,Li2017ConEng,Li2017,Liu2011,lu2016controllability,tian2018controllability,wang2009controllability,Yan2012PRL}.
By detecting the controllability of the underlying networks, one can implement various control tasks to alter systems' states accordingly.
To do this, the associated  energy cost, serving as  a common metric, has to be estimated in advance.

In many large-scale practical dynamical networks, it is a strong constraint to ensure their full controllability.
Moreover, in practical control tasks, it is prohibitively expensive and unnecessary to steer the whole network nodes towards the desired state.
In Ref.~\cite{Gao2014}, authors approximated the minimum number of driver nodes for the target control of complex networks, which shows traditional full control overestimates the number of driver nodes.
That is, choosing a subset of network nodes as target nodes and only controlling these nodes to achieve expected tasks efficiently reduces the number of control inputs.
For energy cost of target control, Ref. \cite{Klickstein2017} presented that it exponentially ascends with the number of target nodes.
Therefore, it requires much more energy to control the entire network.
However, previous analysis depends on the full controllability of the entire network, namely, although we only need to calculate the energy cost to control some of the nodes of the network, other nodes have to be controllable as well \cite{Klickstein2017}.
A systematic analysis of energy cost for achieving the sole target control with the existence of uncontrollable nodes remains elusive.

Here, we consider energy cost for target control.
We present the scaling behavior of both upper and lower bounds of the minimum energy in terms of the given control time.
Furthermore, we revel that with the certain number of target nodes, \text{different} targets can result in hugely different energy cost.
Particularly, for the case of nonsingular adjacency matrix, the corresponding results are more intuitively over the networked topology, where we just need to examine the induced subgraph spanned by the target nodes.

\section{Results}
\begin{figure}[!ht]
\centering
\begin{minipage}[t]{1\textwidth}
\centering
\includegraphics[width=15cm]{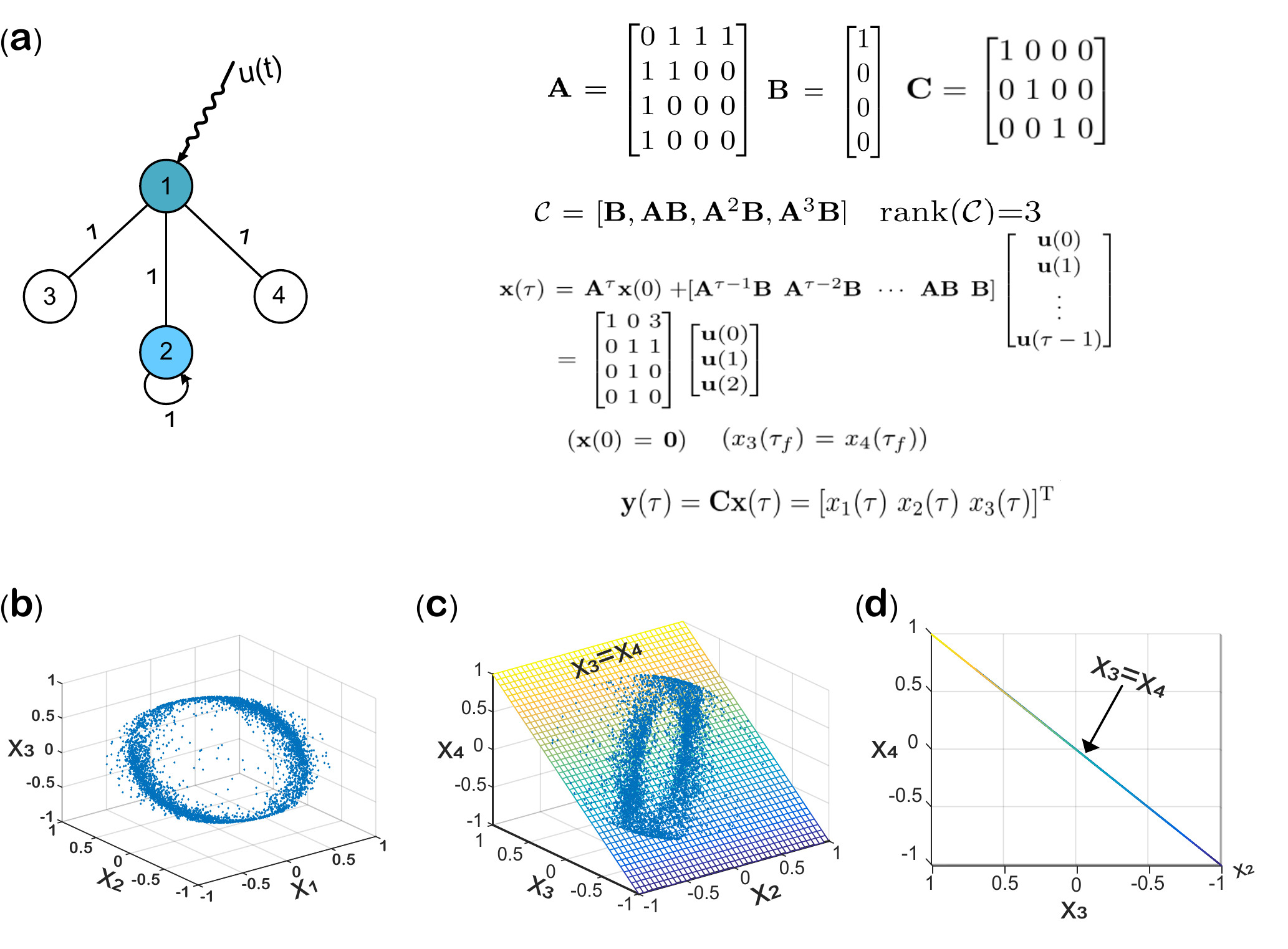}
\end{minipage}
\caption{Illustration of the controllable space as we choose node $1$ as driver node. (a) A network with $4$ nodes, which is not fully controllable. %due to dilation.
We randomly generate $5, 000$ normalized states that the network can be driven from the origin.
In (b), we depict locations of nodes $1, 2, 3$ in  three dimensional space.
And the corresponding locations of nodes $2, 3, 4$ are presented in (c), from which we can see these points are located in the plane $x_3=x_4$.
Panel (d) is the view of (c) from the angle parallel to the plane $x_3=x_4$.
In other words, (d) presents the projection of (c) onto  the plane $x_2=0$,
which indicates the states of nodes $3$ and $4$ are exactly the same. }\label{figure2}
\end{figure}

We consider the canonical linear discrete time-invariant dynamics
\begin{equation}\label{sys1}
\1(\tau+1)=\A\1(\tau)+\B\U(\tau), \quad \tau=0, 1, 2,...
\end{equation}
where $\1(\tau)=(x_1(\tau), x_2(\tau), \dots, x_n(\tau))^{\T}\in\R^{n}$ denotes the state of the entire network with $n$ nodes and $\U(\tau)=(u_1(\tau), u_2(\tau),$ $\dots, u_m(\tau))^{\T} $ $\in\R^{m}$  captures $m$ external input signals at time $\tau$.
$\A\in\R^{n\times n}$ is the adjacency matrix of the network which can represent \text{interactions} between system components.
Here, we consider undirected networks i.e., $\A=\A^{\T}$.
$\B=(b_{ij})\in\R^{n\times m}$ is the input matrix, in which $b_{ij}=1$ means node $i$ is infected by the input $j$ directly, otherwise $b_{ij}=0$.
Nodes that directly receive independent external input signals are called driver nodes.
And, each external input is allowed to directly control one driver node (Fig.~\ref{figure2}(a)).
System (\ref{sys1}) is defined as controllable in $\tau_f\in \N$ steps if there exists an input $\U(\tau)$ to drive the system to the desired state $\1_f=\1(\tau_f)$ from any initial state $\1_0=\1(0)$ at $\tau_f$ steps.
According to the Kalman's controllability criterion \cite{Kalman63}, the rank $r=\text{rank}(\mathcal{C})$ of the controllability matrix
\begin{equation*}\label{contrllability matrix}
\mathcal{C}(\A, \B)=[\B, \A\B, \A^2\B, \cdots,\A^{n-1}\B]
\end{equation*}
represents the dimension of the controllable space.
In other words, for networked system (\ref{sys1}), $r=\text{rank}(\mathcal{C})$ tells that there are $r$ nodes can be controlled towards any desired state in finite time.
For example, for a not fully controllable network shown in fig.~\ref{figure2}(a) with $4$ nodes, the controllable space is three dimensional.
From fig.~\ref{figure2}(c), (d), we can see that the states of nodes $3$ and $4$ are always the same.
The corresponding method determining the controllable subspace was proposed in \cite{hosoe1980determination}.
In what follows, we aim to analyse the corresponding energy cost required to control these controllable nodes.

Firstly, we consider the special case where  $\A$ is nonsingular.
We analyse the target control of this case from the view of structural controllability \cite{Lin1970,Liu2011}.
For a network, if each node has a self dynamics (i.e., $a_{ii}\neq 0$), then the corresponding adjacency matrix is nonsingular.
And for the data-sample system, there are two types of signals:  continuous-time signal and discrete-time signal.
In order to study and design this kind of system, it is necessary to translate the continuous-time state space models into the equivalent discrete-time state space models \cite{Zhenglinear}.
Specifically, for continuous-time systems
\begin{equation*}\label{cts}
\dot{\1}(t)=\mathcal{A} \1(t)+\mathcal{B} \U(t),
\end{equation*}
where $\mathcal{A}$ and $\mathcal{B}$ are the corresponding system matrix and input matrix,
the \text{equivalent} discretized linear system is the system (\ref{sys1}),
where $\1(\tau)=[\1(t)]_{t=\tau\eta},$ $\U(\tau)=[\U(t)]_{t=\tau\eta}, \A=\e^{\mathcal{A} \eta}$ and $\B=(\int^{\eta}_0\e^{\mathcal{A} t}\text{d}t)\mathcal{B}$ with $\eta$ being equal sampling period satisfying Shannon sampling theorem \cite{shannon1998communication}.
Apparently, in the time-discretized system (\ref{sys1}), the system matrix $\A=\e^{\mathcal{A} \eta}$ is nonsingular.
When system (\ref{sys1}) is not fully controllable, i.e., $\text{rank}(\mathcal{C})=r<n$, we denote the controllable node set by $\{{\text{c}_1},$ ${\text{c}_2}, \cdots, {\text{c}_r}\}$ and uncontrollable node set by $\{{\bar{\text{c}}_1},$ ${\bar{\text{c}}_2}, \cdots, {\bar{\text{c}}_{n-r}}\}$.
Next, we permute  the order of the nodes such that the first $r$ nodes are controllable, and leave the rest $n-r$ nodes uncontrollable.
Consequently, the corresponding state variable transformation comes with
\begin{equation*}\label{pt}
\bar{\1}(\tau)=\X\1(\tau),
\end{equation*}
where $\X$ is a permutation matrix with $\X\X^{\T}=\X^{\T}\X=\I$.
Then system (\ref{sys1}) is equivalent to
\begin{equation*}
\bar{\1}(\tau+1)=\X\A\X^{\T}\bar{\1}(\tau)+\X \B\U(\tau),
\end{equation*}
which  can be further decomposed with the following expression
\begin{equation}\label{eq3}
\begin{split}
\begin{bmatrix}
\bar{\1}_{\text{c}}(\tau+1)\\
\bar{\1}_{\bar{{\text{c}}}}(\tau+1)
\end{bmatrix}
=
\begin{bmatrix}
\X_{\text{c}}\A\X^{\T}_{\text{c}}&\,\,\X_{\text{c}}\A\X^{\T}_{\bar{{\text{c}}}}\\
\X_{\bar{{\text{c}}}}\A\X^{\T}_{\text{c}}&\,\,\X_{\bar{{\text{c}}}}\A\X^{\T}_{\bar{{\text{c}}}}
\end{bmatrix}
\begin{bmatrix}
\bar{\1}_{\text{c}}(\tau)\\
\bar{\1}_{\bar{{\text{c}}}}(\tau)
\end{bmatrix}
+
\begin{bmatrix}
\B_{\text{c}}\\
\mathbf{0}
\end{bmatrix}
\U(\tau).
\end{split}
\end{equation}
Therein, $\bar{\1}_{\text{c}}=[x_{{\text{c}}_1}\,\,x_{{\text{c}}_2}\,\, \dots\,\, x_{{\text{c}}_r}]^{\T}$, $\bar{\1}_{\bar{{\text{c}}}}=[x_{\bar{{\text{c}}}_1}\,\, x_{\bar{{\text{c}}}_2}\,\, \dots\,\, x_{\bar{{\text{c}}}_{n-r}}]^{\T}$, $\X_{\text{c}}$ is the first $r$ rows of $\X$, and $\X_{\bar{{\text{c}}}}$ is the remaining $n-r$ rows.
For $\X_{\text{c}}\A\X^{\T}_{\bar{{\text{c}}}}$, the element of the $i$th row and $j$th column is $\X_{\text{c}}\A\X^{\T}_{\bar{{\text{c}}}}(i, j)=a_{{\text{c}}_i, \bar{{\text{c}}}_j}.$
In other words, we need to judge whether there exists a link between the controllable node ${\text{c}}_i$ and the uncontrollable node $\bar{{\text{c}}}_j$.
A network is structural controllable, if the following conditions are satisfied \cite{Lin1970,Liu2011}: 1) each node is accessible from external inputs; 2) there is no dilation.
In addition, if matrix $(\A~ \B)$ has full rank, then the network $(\A, \B)$ has no dilation.
Here $\A$ is nonsingular which leads to the absence of dilation  in the network $(\A, \B)$.
If there exists a direct link between nodes ${\text{c}}_i$ and $\bar{{\text{c}}}_j$, then the node $\bar{{\text{c}}}_j$ is controllable,
since $\bar{{\text{c}}}_j$ can receive control signal from the controllable ${\text{c}}_i$.
Therefore, there is no direct link between ${\text{c}}_i$ and  $\bar{{\text{c}}}_j$, i.e., $a_{{\text{c}}_i, \bar{{\text{c}}}_j}=0$.
By that, system (\ref{eq3}) is
\begin{equation*}\label{eq030704}
\begin{bmatrix}
\bar{\1}_{\text{c}}(\tau+1)\\
\bar{\1}_{\bar{{\text{c}}}}(\tau+1)
\end{bmatrix}
=
\begin{bmatrix}
\X_{\text{c}}\A\X^{\T}_{\text{c}}&\,\,\0\\
\0&\,\,\X_{\bar{{\text{c}}}}\A\X^{\T}_{\bar{{\text{c}}}}
\end{bmatrix}
\begin{bmatrix}
\bar{\1}_{\text{c}}(\tau)\\
\bar{\1}_{\bar{{\text{c}}}}(\tau)
\end{bmatrix}
+
\begin{bmatrix}
\B_{\text{c}}\\
\0
\end{bmatrix}
\U(\tau).
\end{equation*}
And the controllable subsystem is
\begin{equation*}\label{eq4}
\bar{\1}_{\text{c}}(\tau+1)=\X_{\text{c}}\A\X^{\T}_{\text{c}}\bar{\1}_{\text{c}}(\tau)+\B_{\text{c}}\U(\tau).
\end{equation*}
By letting $\bar{\A}_{\text{c}}=\X_{\text{c}}\A\X^{\T}_{\text{c}}$, it is clear that all eigenvalues of $\bar{\A}_{\text{c}}$ belong to the set of eigenvalues of $\A$.
And the graph of the adjacency matrix $\bar{\A}_{\text{c}}$ is a subgraph of original graph, spanned by nodes ${{\text{c}}_1}, {{\text{c}}_2}, \dots, {{\text{c}}_r}$.
Note that, a network is structural controllable can be a precondition of state controllable.
And for a structural controllable network, almost all weights can guarantee the controllability of states \cite{Lin1970}.
In this part, we do not consider the rare scenarios where the network is structural controllable but not state controllable.
In this case, we obtain the energy cost for controlling a controllable network in table \ref{table2} of Section~\ref{Energycom}.

For  the general case, we can equivalently regard target control as the output controllability of a system.
Namely, consider a part of the system (\ref{sys1})
\begin{equation}\label{sys2}
\begin{cases}
\1(\tau+1)=\A\1(\tau)+\B\U(\tau)\\
\y(\tau)=\C\1(\tau)
\end{cases}
\end{equation}
where $\C=[\I_{\text{c}_1}^{\T}~\I_{\text{c}_2}^{\T}~\cdots~\I_{\text{c}_r}^{\T}]^{\T}\in \R^{r\times n}$ is the output matrix with $\I_{\text{c}_i}$ being the $i$th row of identity matrix.
$\y(\tau)%=[y_1(tau)~y_2(t)~\cdots~y_{r}(t)]^{\T}
=[x_{{\text{c}}_1}(\tau)\,\,x_{{\text{c}}_2}(\tau)\,\, \dots\,\, x_{{\text{c}}_r}(\tau)]^{\T}$ collects the states of target nodes.
In fig.~\ref{figure2}, we make a brief explanation on target control.
System (\ref{sys2}) is called output controllable  if and only if the output controllability matrix satisfies
\begin{equation}
\text{rank}~ \mathcal{C}(\A, \B, \C)=\text{rank}~[\C\B, \C\A\B, \C\A^2\B, \cdots,\C\A^{n-1}\B]=r.
\end{equation}
Note that partial controllability of the system (\ref{sys1}) is equivalent to output controllability of the system (\ref{sys2}) \cite{lin1991output}.
Without loss of generality, we denote controllable nodes by $1, 2, \dots, r$, i.e., $\C$ is chosen from the first $r$ rows of an identity matrix with size $n$.
To analyse the energy cost for target control, we employ the conventional definition of the following input control energy
\begin{equation}\label{ET}
E(\tau_f)=\frac{1}{2}\sum^{\tau_f-1}_{\tau=0}\U^{\T}(\tau)\U(\tau).
\end{equation}
By minimizing the energy cost $E(\tau_f)$, one can employ optimal energy control theory \cite{OptimalBooLewis} to derive the optimal control input (see section \ref{outcon})
\begin{equation}\label{opU}
\U^*(\tau)=\B^{\T}(\A^{\T})^{\tau_f-\tau-1}\C^{\T}(\C\W\C^\T)^{-1}(\y_{f}-\C\A^{\tau_f}\1_0),
\end{equation}
where $\y_f=\y(\tau_f)$ and $\W$ is Gramian matrix of the system (\ref{sys1}) with
$
\W=\sum^{\tau_f-1}_{i=0}\A^{\tau_f-i-1}\B$ $\B^{\T}(\A^{\T})^{\tau_f-i-1}.
$
Accordingly, the minimum energy cost is
\begin{equation}\label{minE}
E(\tau_f)=(\y_f-\C\A^{\tau_f}\1_0)^{\T}(\C\W\C^{\T})^{-1}(\y_f-\C\A^{\tau_f}\1_0).
\end{equation}
Assuming $\1_0=\mathbf{0}$ and denoting $\W_{\text{C}}=\C\W\C^{\T}$, we further have
\begin{equation}
E(\tau_f)=\y_f^{\T}\W_{\text{C}}^{-1}\y_f.
\end{equation}
Intuitively, $\W_{\text{C}}$ consists of the first $r$ rows and $r$ columns of the matrix $\W$.
Note that if system (\ref{sys2}) is output controllable, matrix $\W_{\text{C}}$ is invertible.
By normalizing the control distance  $\|\y_f\|=1$, we have
\begin{equation}\label{ineqe}
\frac{1}{\lambda_{\max}(\W_{\text{C}})} \leq E(\tau_f)\leq \frac{1}{\lambda_{\min}(\W_{\text{C}})},
\end{equation}
where $\lambda_{\max} (\lambda_{\min})$ is the maximum (minimum) eigenvalue of $\W_{\text{C}}$.
From Eq.~(\ref{ineqe}), the kernel problem is to obtain the minimum and maximum eigenvalues of $\W_{\text{C}}$.

\begin{figure}[ht]
\centering
\begin{minipage}[t]{1\textwidth}
\centering
\includegraphics[width=15cm]{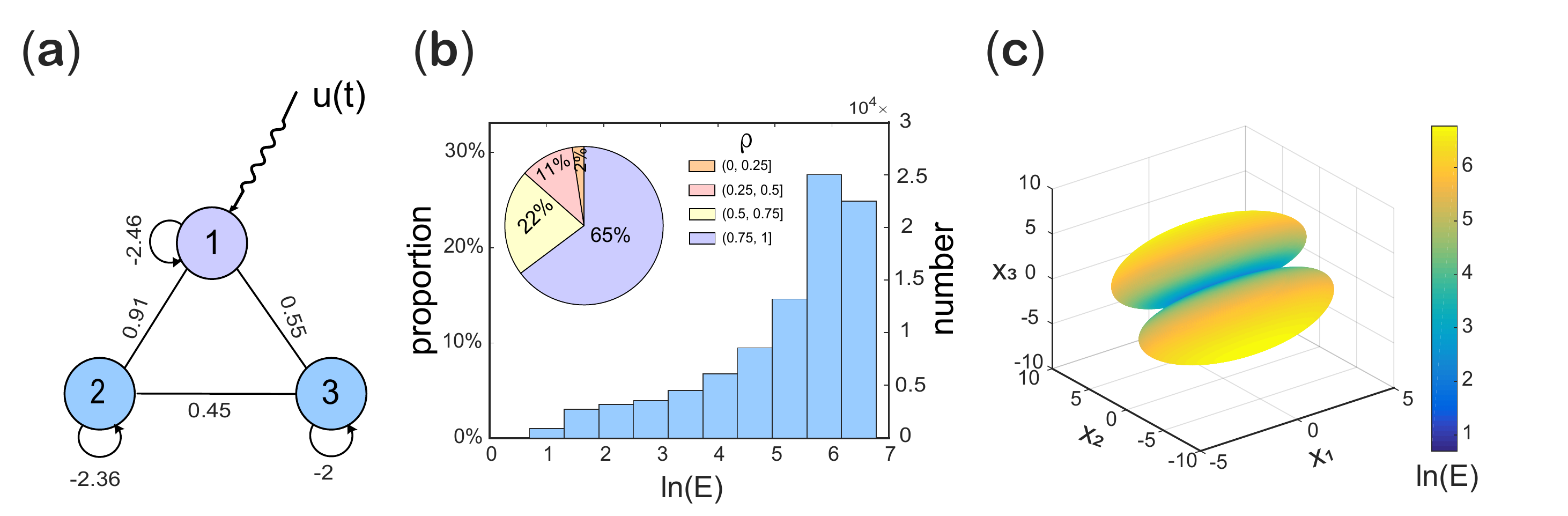}
\end{minipage}
\caption{Energy distribution for a three dimensional fully controllable network.
(a) For the fully connected three dimensional controllable network, the interaction strengths are given alongside each link, and node $1$ is chosen as driver node.
In (b), we  make the statistics of energy distribution for driving a fully controllable network from the origin to unit ball by traversing about $90,000$ points from the surface uniformly.
For these $90,000$ data, we find the maximum and the minimum values and divide them into 10 equal intervals.
The bar diagram counts the number and percentage of the data in each subinterval.
In addition, we calculate the ratios of each data to the maximum value, denoted by $\rho$.
Therefore, we count the probability of the ratio falling into the interval $(0, 0.25], (0.25, 0.5], (0.5, 0.75]$ and $(0.75, 1]$, as shown in the pie.
Accordingly in (c), we depict these $90,000$ points in three dimensional coordinates by taking  logarithm of these energy ln$(E)$ as the distance from the origin to these points.}\label{figure1}
\end{figure}

In previous studies,  researchers focused on the lower bound of energy cost, where  the trace of the corresponding Gramian matrix used to approximate the maximum eigenvalue \cite{Li2017,Pasqualetti2014Controllability,Yan2012PRL}.
Indeed, given that the trace of a matrix is equal to the sum of the corresponding eigenvalues, it is frequently employed to approximate the maximum eigenvalue of controllability Gramian matrix, where most eigenvalues are relatively small.
It further reflects the upper bound of the energy cost dominates.
Therefore, it is meaningful to acquire the corresponding upper bound of energy for achieving control goals.
We take the three-dimensional fully controllable network as an example.
As shown in Eq.~(\ref{eq5}), matrix $\W^{-1}$ of $\alpha^{\T}\W^{-1}\alpha$ is invertible with three positive eigenvalues $\mu_1, \mu_2, \mu_3$ and the corresponding three linearly independent and orthogonal normalized eigenvectors $\alpha_1, \alpha_2, \alpha_3$.
Here, we assume $\|\alpha\|=1$ and $\alpha$ is a linear combination of $\alpha_1, \alpha_2, \alpha_3$ as $\alpha=a_1\alpha_1+a_2\alpha_2+a_3\alpha_3$.
Therefore, we have $E=\mu_1a_1^2+\mu_2a_2^2+\mu_3a_3^2$ with $a_1^2+a_2^2+a_3^2=1$.
After introducing new variables $x=\sqrt{\mu_1}a_1$, $y=\sqrt{\mu_2}a_2$ and $z=\sqrt{\mu_3}a_3$ with $a_1=\sin\theta\cos\phi$, $a_2=\sin\theta\sin\phi$ and $a_3=\cos\theta$,
we have $x^2+y^2+z^2=E(\alpha)=d^2$, where $d$ is the distance between the origin and the point $(x, y, z)$.
At a given control distance, we show the energy required to reach final states on the unit sphere in  fig.~\ref{figure1}(c).
And from fig.~\ref{figure1}, it is clear that a large proportion ($\sim 65\%$) of the final states requires the largest amount ($>75\%$) of energy to accomplish control.

In the sequel, we analyse the problem from the point of view of system decomposition.
Since system (\ref{sys1}) is not fully controllable, the system can be decomposed into two parts: controllable part and uncontrollable part,
by introducing variable transformation
\begin{equation}
\overline{\1}(\tau)=\2\1(\tau),
\end{equation}
where $\2$ is  an orthogonal matrix.  The first $r$ columns of $\2^{\T}$ are constructed by the orthonormal basis of column space of $\mathcal{C}$ (via Gram-Schmidt orthogonalization), and the rest columns are constructed by $n-r$ column vectors orthogonal to existing $r$ columns.
Therefore, system (\ref{sys1})  translates into
\begin{equation}\label{eq1}
\overline{\1}(\tau+1)=\overline{\A}\overline{\1}(\tau)+\overline{\B}\U(\tau)
\end{equation}
where $\overline{\A}=\2\A\2^{\T}$ and $\overline{\B}=\2\B$.
According to controllability decomposition theory, the specific form of (\ref{eq1}) is
\begin{equation}\notag
\begin{bmatrix}
\overline{\1}_{\text{c}}(\tau+1)\\
\overline{\1}_{\text{nc}}(\tau+1)
\end{bmatrix}
=
\begin{bmatrix}
\A_{\text{c}}&\mathbf{0}\\
\mathbf{0}&\A_{\text{nc}}
\end{bmatrix}
\begin{bmatrix}
\overline{\1}_{\text{c}}(\tau)\\
\overline{\1}_{\text{nc}}(\tau)
\end{bmatrix}
+
\begin{bmatrix}
\B_{\text{c}}\\
\mathbf{0}
\end{bmatrix}
\U(\tau),
\end{equation}
where $\overline{\1}_{\text{c}}(\tau)\in \R^{r}$, $\A_{\text{c}}\in\R^{r\times r}$, and $\B_{\text{c}}\in \R^{r\times m}$.
Therein, the dynamics of the controllable part is
\begin{equation}\label{eq2}
\overline{\1}_{\text{c}}(\tau+1)=\A_{\text{c}}\overline{\1}_{\text{c}}(\tau)+\B_{\text{c}}\U(\tau).
\end{equation}

For example, for a not fully controllable network as shown in fig.~\ref{figure2}, one can preform controllable decomposition.
Firstly, it needs to obtain maximal linearly independent group of  column vectors of $\mathcal{C}$ as
$
\protect\begin{bmatrix}
1\protect\\
0\protect\\
0\protect\\
0
\protect\end{bmatrix},
\protect\begin{bmatrix}
0\protect\\
1\protect\\
1\protect\\
1
\protect\end{bmatrix},
\protect\begin{bmatrix}
3\protect\\
1\protect\\
0\protect\\
0
\protect\end{bmatrix}.
$
By performing Gram-Schmidt orthogonalization on the above vectors and extending to the whole $4$ dimensional linear space, we have
$\2^{\T}=\protect\begin{bmatrix}
1&0&0&0\protect\\
0&\frac{1}{\sqrt{3}}&\frac{2}{\sqrt{6}}&0\protect\\
0&\frac{1}{\sqrt{3}}&\frac{-1}{\sqrt{6}}&\frac{1}{\sqrt{2}}\protect\\
0&\frac{1}{\sqrt{3}}&\frac{-1}{\sqrt{6}}&\frac{-1}{\sqrt{2}}
\protect\end{bmatrix}.$
Furthermore,
$$\overline{\A}=\2\A\2^{\T}=
\protect\begin{bmatrix}
0&\sqrt{3}&0&0\protect\\
\sqrt{3}&\frac{1}{3}&\frac{2}{\sqrt{18}}&0\protect\\
0&\frac{2}{\sqrt{18}}&\frac{2}{3}&0\protect\\
0&0&0&0
\protect\end{bmatrix}
\quad \text{with} \quad \A_{\text{c}}=\protect\begin{bmatrix}
0&\sqrt{3}&0\protect\\
\sqrt{3}&\frac{1}{3}&\frac{2}{\sqrt{18}}\protect\\
0&\frac{2}{\sqrt{18}}&\frac{2}{3}
\protect\end{bmatrix}$$
and
$\B_{\text{c}}=[1 ~ 0~0]^{\T}.$

For system (\ref{eq2}), the corresponding Gramian matrix  is
\begin{equation}
\mathcal{W}=\sum^{\tau_f-1}_{\tau=0}\A^\tau_{\text{c}}\B_{\text{c}}\B_{\text{c}}^{\T}\A_{\text{c}}^\tau,
\end{equation}
which is invertible.
Substituting $\overline{\A}=\2\A\2^{\T}$ and $\overline{\B}=\2\B$ into $\W_{\text{C}}$, we have
\begin{equation}\label{WCT}
\W_{\text{C}}=\2_1^{\T}\mathcal{W}\2_1,
\end{equation}
where $\2=
\begin{bmatrix}
\2_1& \2_3\\
\2_2& \2_4
\end{bmatrix}$
and $\2_1\in \R^{r\times r}$.
It needs to be emphasized that the ultimate goal is to derive the minimum and maximum eigenvalues of $\W_{\text{C}}$.
For that, based on an effective approach to approximate the minimum and maximum eigenvalues of positive definite matrix $\M$ \cite{lam2011estimates}, we can derive the corresponding  upper and lower bounds as
\begin{equation}\label{upE}
\overline{E}\approx f(\underline{\alpha}, \underline{\beta})
\end{equation}
\begin{equation}\label{lowE}
\underline{E}\approx \frac{1}{f(\overline{\alpha}, \overline{\beta})}
\end{equation}
where
$f(\alpha, \beta)=\sqrt{\frac{\alpha}{n}+\sqrt{\frac{n-1}{n}(\beta-\frac{\alpha^2}{n})}}$,
$
\overline{\alpha}=\text{trace}(\M^2),
\overline{\beta}=\text{trace}(\M^4),
\underline{\alpha}$= \text{trace} $((\M^{-1})^2),
$
and
$
\underline{\beta}=\text{trace}((\M^{-1})^4).
$

\begin{figure}[http]
\centering
\includegraphics[width=15cm]{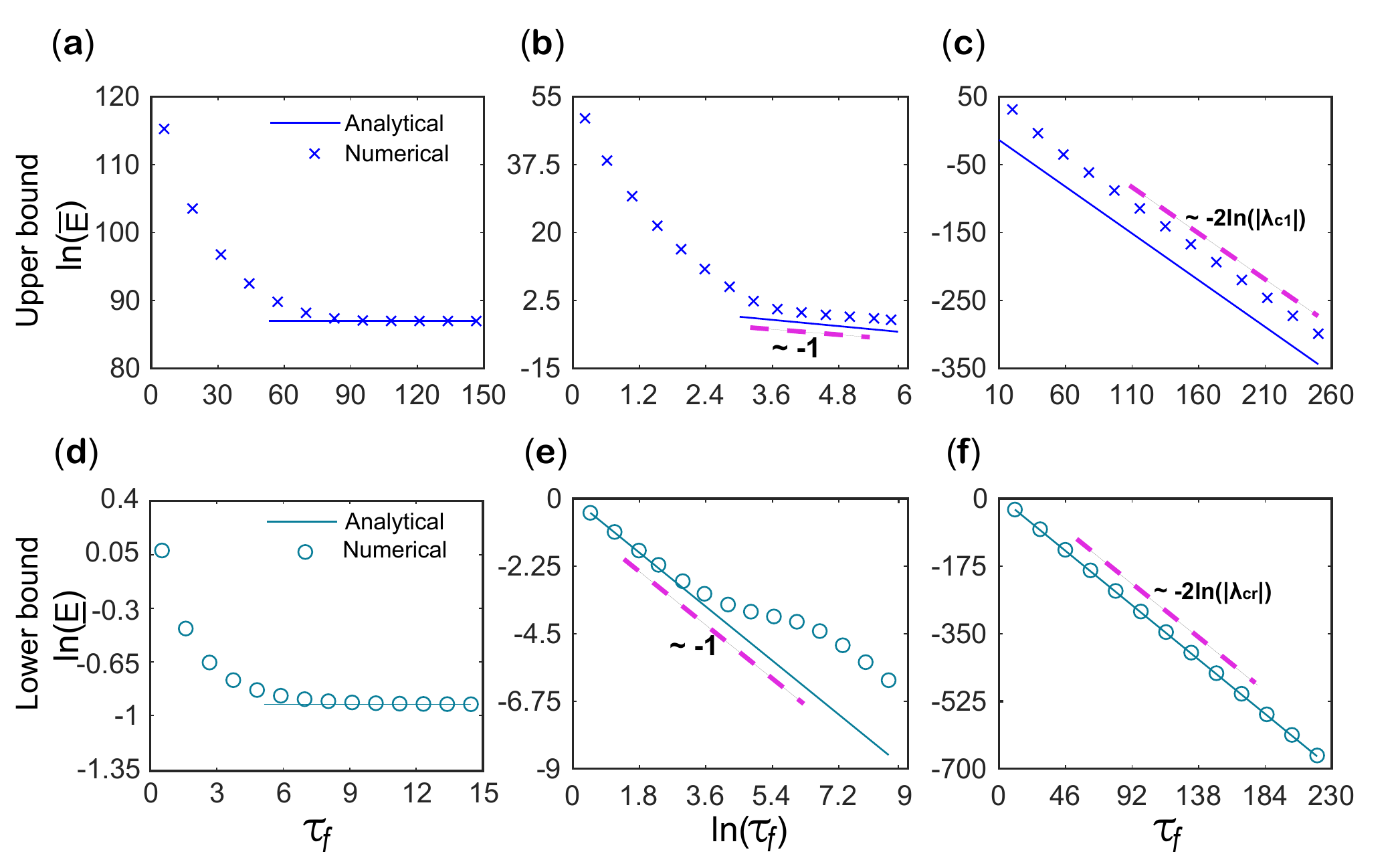}
\caption{The upper and lower bounds of control energy for target control.
In order to show the theoretical results in table~\ref{table1}, we generate random networks with different types of $\A_{\text{c}}$.
For each pair of nodes, we add an edge with the probability $0.1$ \cite{erdHos1960evolution}.
To systematically show all cases given in table~\ref{table1}, in (a) and (d), here we set link weight $a_{ij}$ from the interval $[0, 0.05]$ uniformly.
Analogously, in (b), (c) and (f), the interval is $[-1, 0]$, and in (e), the interval is $[0, 0.1]$.
To intuitively judge the eigenvalues of $\A_{\text{c}}$, we add self-loops and set the corresponding weight as $a+s_i$ with $s_i=-\sum^n_{j=1}a_{ij}$.
In (a) and (d), by setting  $a=0.8$, we can derive that all eigenvalues of $\A_{\text{c}}$ are located in $[0.6478, 0.8]$, which leads to $|\lambda_{\text{c}1}|<1$ and $|\lambda_{\text{c}r}|<1$.
In (b), (e) and (f), we set $a=1$ such that all eigenvalues of $\A_{\text{c}}$ of (b) and (f) are located in $[1, 4.4048]$ with $|\lambda_{\text{c}1}|=1$ and $|\lambda_{\text{c}r}|>1$ and all eigenvalues of $\A_{\text{c}}$ of (e) are located in $[0.5805, 1]$ with $|\lambda_{\text{c}r}|=1$.
In (c), all eigenvalues of $\A_{\text{c}}$ are located in $[2, 6.5169]$ by selecting $a=2$, which leads to $|\lambda_{\text{c}1}|>1$.
Here, we adopt the random network with size $n=20$, and the dimension of controllable space of the network is $14$ by selecting one driver node.
In each panel, we can see that the generated pattern of analytical derivation almost overlaps that of numerical calculations.
}\label{figure7}
\end{figure}

In Sections \ref{Energycom} and \ref{energyuncom}, we make adequate analysis on calculating parameters  $\overline{\alpha}$, $\overline{\beta}$, $\underline{\alpha}$ and $\underline{\beta}$ of $\W_{\text{C}}$ to obtain approximations of the corresponding maximum and minimum eigenvalues.
Table~\ref{table1} presents the upper and lower bounds of the minimum energy cost and the corresponding numerical verification is shown in fig.~\ref{figure7}.

\begin{table}[!http]
\centering\caption{Lower and upper bounds of the minimum energy.
Here, $|\lambda_{\text{c}1}|$ and $|\lambda_{\text{c}r}|$ are the minimum and the maximum absolute values of eigenvalues of $\A_{\text{c}}$, respectively, i.e., $|\lambda_{\text{c}1}|=\min\{|\lambda_{\text{c}i}|~|\lambda_{\text{c}i}\in \lambda(\A_{\text{c}})\}$ and $|\lambda_{\text{c}r}|=\max\{|\lambda_{\text{c}i}|~|\lambda_{\text{c}i}\in \lambda(\A_{\text{c}})\}$.
When $|\lambda_{\text{c}r}|>1$, $\underline{E}$ decreases for time $\tau_f$ with power law $\lambda_{\text{c}r}^{2-2\tau_f}$.
When $|\lambda_{\text{c}r}|=1$, $\underline{E}\sim \tau_f^{-1}$ holds for any number of driver nodes.
When $|\lambda_{\text{c}r}|=1$, $\underline{E}$ approaches a constant irrespective of $\tau_f$.
For one driver node, the constant is given by Eq.(\ref{lowE}) with (\ref{Aeq11})(\ref{Aeq12}) and for any number of driver nodes, the constant is given by Eq.(\ref{lowE}) with (\ref{Aeq15})(\ref{Aeq16}).
Analogously, for the upper bound, when  $|\lambda_{\text{c}1}|>1$, $\overline{E}\sim \lambda_{\text{c}1}^{2-2\tau_f}$ holds and when $|\lambda_{\text{c}1}|=1$, $\overline{E}\sim \tau_f^{-1}$ holds. When $|\lambda_{\text{c}1}|<1$, $\overline{E}$ approaches a constant irrespective of $\tau_f$ as Eq.(\ref{upE}) with (\ref{Aeq13})(\ref{Aeq14}) for one driver node.
}
\fontsize{8}{15}\selectfont
\begin{tabular}{cc|c|c}
\toprule[2pt]
\multicolumn{2}{c|}{Number of driver nodes}&$1$&$m (m\leq r)$\\
\toprule[1pt]
%\hline
\multirow{3}{*}{Lower bound $\underline{E}$}&$|\lambda_{\text{c}r}|<1$&Eq.(\ref{lowE}) with (\ref{Aeq11})(\ref{Aeq12})&Eq.(\ref{lowE}) with (\ref{Aeq15})(\ref{Aeq16})\\
%\hline
&$|\lambda_{\text{c}r}|=1$&$ \sim \tau_f^{-1}$& $\sim \tau_f^{-1}$\\
%\hline
&$|\lambda_{\text{c}r}|>1$&$\sim \lambda_{\text{c}r}^{2-2\tau_f}$&$\sim \lambda_{\text{c}r}^{2-2\tau_f}$\\
\toprule[1pt]
\multirow{3}{*}{Upper bound $\overline{E}$}&$|\lambda_{\text{c}1}|<1$&Eq.(\ref{upE}) with (\ref{Aeq13})(\ref{Aeq14})&$C_1$\\
%\hline
&$|\lambda_{\text{c}1}|=1$&$ \sim \tau_f^{-1}$& $\sim \tau_f^{-1}$\\
%\hline
&$|\lambda_{\text{c}1}|>1$&$\sim \lambda_{\text{c}1}^{2-2\tau_f}$&$\sim \lambda_{\text{c}1}^{2-2\tau_f}$\\
\toprule[2pt]
\end{tabular}
\label{table1}
\end{table}

\begin{figure}[ht]
\centering
\includegraphics[width=14cm]{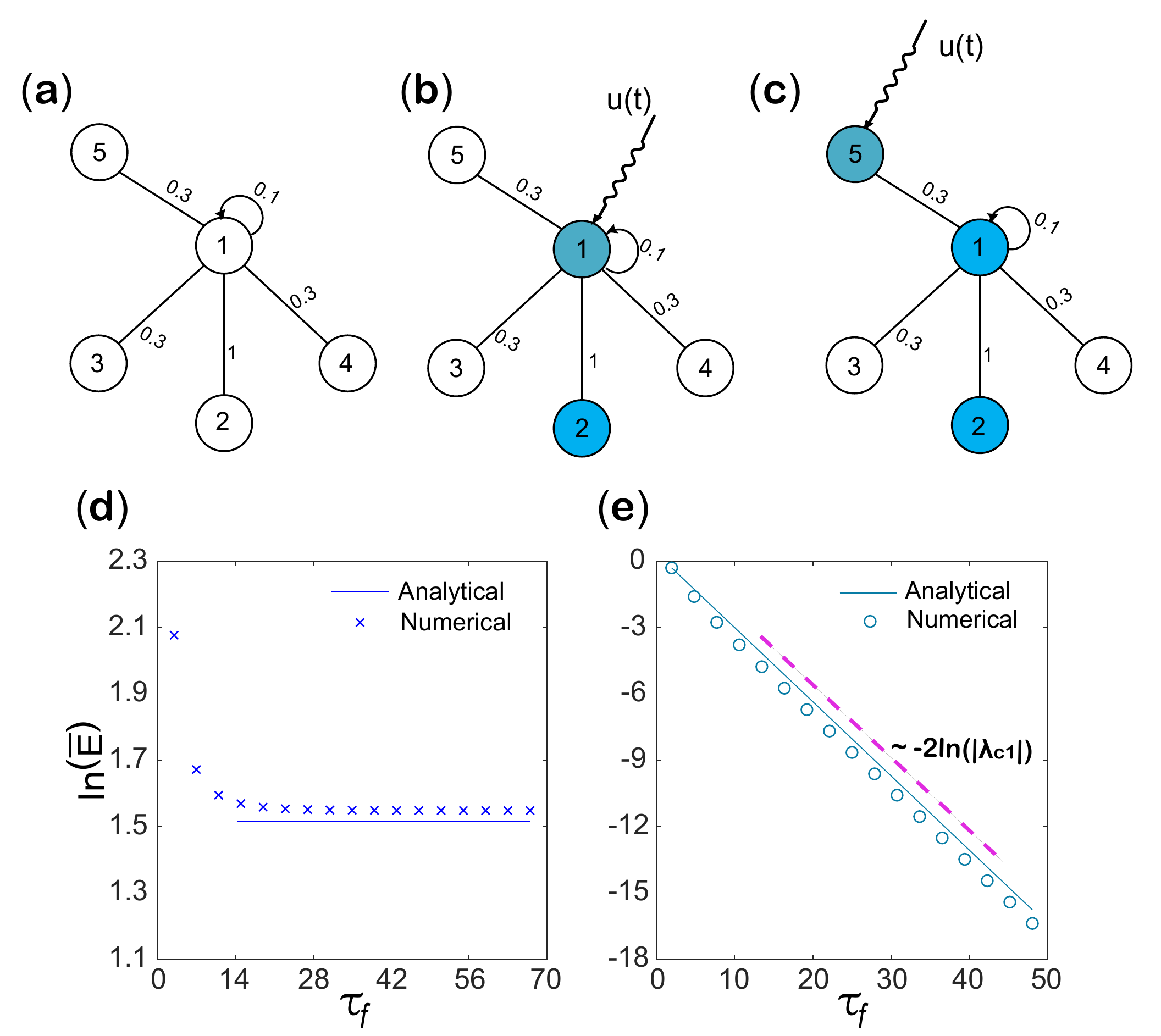}
\caption{Target control under different choices of driver nodes.
(a) A network with $5$ nodes.
In (b) and (c), we assign node $1$ and node $5$ as driver node, respectively.
When node $1$ is driver node, only node $1$ and node $2$ are controllable; when node $5$ is driver node, nodes $1, 2, 5$ are controllable.
For (b), the minimum eigenvalue $\lambda_{\text{c}1}$ of $\A_{\text{c}}$ satisfies $|\lambda_{\text{c}1}|<1$.
Therefore, the corresponding upper bound of energy approaches a constant, depicted in (d) (see theory shown in table~\ref{table1}).
Analogously, for (c), the minimum eigenvalue $\lambda_{\text{c}1}$ of $\A_{\text{c}}$ satisfies $|\lambda_{\text{c}1}|>1$.
Therefore, the corresponding upper bound of energy exponentially decreases with $\tau_f$, as depicted in (e).
}\label{figure4}
\end{figure}

From energy scaling in terms of the control time for target control, we find that controlling different target nodes corresponds to different energy scalings.
Therefore, in the given control task to control a given number of nodes, one can achieve minimum energy control by choosing appropriate driver nodes.
In fig.~\ref{figure4}(a), it is clear that the network is not controllable due to dilation.
When \text{nodes} $1$ and $5$ are assigned as driver node separately, set of the controllable nodes is different.
When node $1$ is  the driver node,  nodes $1, 2$ are controllable, and the minimum absolute value of eigenvalues  of the corresponding Gramian matrix of controllable part is $0$.
When  node $5$ is  the driver node,  nodes $1, 2, 5$ are controllable, and the minimum absolute value of  eigenvalues of the corresponding Gramian matrix of controllable part is $1.07$.
Therefore, the upper bounds of energy cost for achieving target control are different, and the corresponding energy scaling behaviors are different, which are depicted in fig.~\ref{figure4}(d), (e).
For a given set of target nodes, one can employ a greedy algorithm proposed in \cite{Gao2014} to find an approximately minimum set of driver nodes for target control.

\begin{figure}[ht]
\centering
\begin{minipage}[t]{1\textwidth}
\centering
\includegraphics[width=15cm]{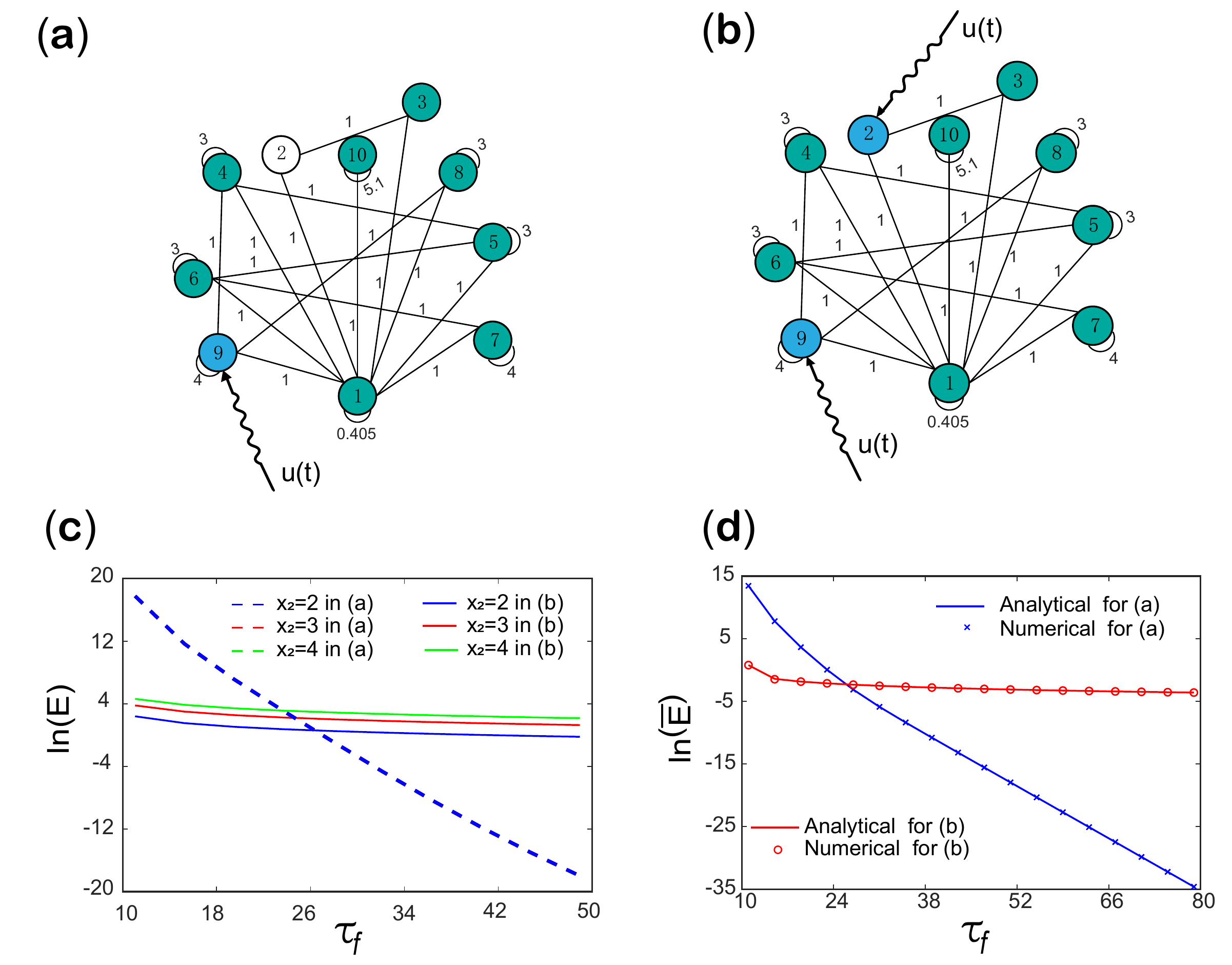}
\end{minipage}
\caption{Comparison of energy cost for achieving target control and full control.
In (a), we control the network by adding an input to node $9$ directly (blue), under which there are $9$ nodes being controllable (green) except the node $2$ (white).
In (b), we directly control nodes $2$ and $9$, by which the network is fully controllable.
In (c), for $\1_f=[x_1(\tau_f)~ x_2(\tau_f)~\cdots~ x_{10}(\tau_f)]^{\T}$ in networks (a) and (b), by letting $x_i(\tau_f)=1, i=1, 3, 4,\dots,10$ and $x_2(\tau_f)=2, 3, 4$, we calculate the corresponding energy cost.
And upper bounds of the corresponding minimum energy cost are presented in (d).
It is clear that energy cost of achieving target control  is much less than achieving full control, when $\tau_f$ is large.
}\label{figure6}
\end{figure}

In addition, we find that we can save huge amount of energy cost, even multiple orders of magnitude less by achieving target control compared with traditional full controllability.
For example, in fig.~\ref{figure6}, we make a comparison on energy cost between controlling partial network and controlling the entire network.
When achieving target control by one external input, there are $9$ controllable nodes as shown in fig.~\ref{figure6}(a) (node $2$ is uncontrollable).
When applying two external inputs to the network as shown in fig.~\ref{figure6}(b), the network is fully controllable.
From fig.~\ref{figure6}(c), we find that it expends less energy cost to achieve target control with controlling $90\%$ nodes, compared with controlling the entire network when $\tau_f$ is large.
When $\tau_f$ is small, it expends less energy cost to achieve full control compared with achieving target control, which is reasonable due to the effect of more input signals for full control and consistent with the result in Ref. \cite{Pasqualetti2014Controllability}.

\section*{Discussion}

In this paper, we investigate energy cost for achieving target control, i.e., controlling a part of complex networks towards the desired state.
In practical control tasks, sometimes it is unnecessary and inadvisable to control the entire network \cite{galbiati2013power,Gao2014}.
Target control relaxes the requirement from full controllability and avoids excessive waste of energy cost \cite{Klickstein2017}.
With respect to issues of energy cost for controlling complex networks,
we give a framework  for estimating the exact upper and the lower bounds of control energy.
The method we applied is effective and can be widely used to analyse the traditional full controllability \cite{Li2018The} as well.

The nonlinear dynamics is admitted to properly describe practical complex systems \cite{Cornelius2013NatCommun,Khalil2002a}.
However, the corresponding operation of recognizing empirical parameterizations is challenging.
A common alternative to process the nonlinearity is to investigate the linearized version.
Indeed, in many aspects, linear dynamics is adequate to approximate and explore nonlinear dynamics in permissible local regions.
Specifically, the ability of detecting the controllability of linearized dynamics has been validated to guarantee the controllability of the original nonlinear \text{systems} along the trajectory in corresponding regions \cite{LinearizationBook}.
Nevertheless, further analysis on nonlinear dynamics still remains a promising direction for general networked \text{systems} in future work.

In  a real complex network, temporal networks are universal, that is, the interactions between nodes are dynamic \cite{Holme2012,Li2017,NaokiBook,Posfai2014NJP}.
As a preliminary exploration, we consider the framework for static networks, which can be extended to explore the corresponding energy cost for controlling  temporal networks.
Controlling a part of temporal networks, one can utilize output controllability to formulate specific problems.
Furthermore, one can take advantage of the effective Gramian matrix given in Ref.~\cite{Li2017} and adopt controllable decomposition to solve the problem.
In order to achieve target control of temporal networks with less energy cost,  one can apply independent path theorem proposed in \cite{Posfai2014NJP} to determine driver nodes according to the framework for temporal networks.

\section*{Acknowledgement}
This work is supported by the National Natural Science
Foundation of China (NSFC) under Grants No. 61751301
and No. 61533001. A.L. acknowledges the support from the
Human Frontier Science Program (HFSP) Postdoctoral Fellowship
(Grant No. LT000696/2018-C) and Foster Lab at Oxford.

%\bibliographystyle{ws-acs}

%\bibliography{bibliography}

\appendix

\section{Optimal Control Energy Theory}

\subsection{Fully controllable systems}

In this subsection, we assume system (\ref{sys1}) is fully controllable.
According to the optimal control energy theory, we have the Hamilton function
\begin{equation}\label{hm}
\h(\tau)=\frac{1}{2}\U(\tau)^{\T}\U(\tau)+\lambda(\tau+1)^{\T}(\A\1(\tau)+\B\U(\tau)).
\end{equation}
The state variable $\1(\tau)$ and $\lambda(\tau)$ satisfy
\begin{equation*}\label{ze1}
\1(\tau+1)=\frac{\partial \h(\tau)}{\partial\lambda(\tau+1)}=\A\1(\tau)+\B\U(\tau),
\end{equation*}
and
\begin{equation}\label{ze2}
\lambda(\tau)=\frac{\partial \h(\tau)}{\partial\1(\tau)}=\A^{\T}\lambda(\tau+1).
\end{equation}
And the optimal input satisfies
\begin{equation}\label{fopu}
0=\frac{\partial \h(\tau)}{\partial\U(\tau)}.
\end{equation}
By solving Eq.~(\ref{fopu}), the optimal input is
\begin{equation}\label{opu1}
\U(\tau)^*=-\B^{\T}\lambda(\tau+1).
\end{equation}
Furthermore, from iterative equation~(\ref{ze2}), we can derive
\begin{equation}\label{Aeq1}
\lambda(0)=(\A^{\T})^\tau\lambda(\tau) \quad \rightarrow\quad  \lambda(0)=(\A^{\T})^{\tau_f}\lambda(\tau_f) \quad \rightarrow\quad  \lambda(\tau)=(\A^{\T})^{\tau_f-\tau}\lambda(\tau_f).
\end{equation}
In addition, the solution of Eq.~(\ref{sys1}) is
\begin{equation}\label{Aeq2}
\1(\tau)=\A^\tau\1_0+\sum^{\tau-1}_{i=0}\A^{\tau-i-1}\B\U(i).
\end{equation}
Substituting Eqs.~(\ref{opu1}) and (\ref{Aeq1}) into (\ref{Aeq2}), we have
\begin{align}
\1(\tau)
&=\A^\tau\1_0-\sum^{\tau-1}_{i=0}\A^{\tau-i-1}\B\B^{\T}(\A^{\T})^{\tau-i-1}\lambda(\tau_f).\label{Aeq3}
\end{align}
By letting $\tau=\tau_f$ in Eq.~(\ref{Aeq3}),  we have
$
\1(\tau_f)=\A^{\tau_f}\1_0-\W\lambda(\tau_f)
$
and
\begin{equation}\label{Aeq4}
\lambda(\tau_f)=-\W^{-1}(\1(\tau_f)-\A^{\tau_f}\1_0)
\end{equation}
where
$
\W=\sum^{\tau_f-1}_{i=0}\A^{\tau_f-i-1}\B\B^{\T}(\A^{\T})^{\tau_f-i-1}
$
is the controllability Gramian matrix of system (\ref{sys1}).
Substituting Eqs.~(\ref{Aeq1}) and (\ref{Aeq4}) into (\ref{opu1}), we derive the optimal input
\begin{equation*}\label{opu}
\U^*(\tau)=\B^{\T}(\A^{\T})^{\tau_f-\tau-1}\W^{-1}(\1(\tau_f)-\A^{\tau_f}\1_0),
\end{equation*}
and the minimum energy cost
\begin{align}
E_{\min}&=\frac{1}{2}\sum^{\tau_f-1}_{\tau=0}\U^*(\tau)^{\T}\U^*(\tau)=\frac{1}{2}\alpha^{\T}\W^{-1}\alpha\label{eq5}
\end{align}
with $\alpha=\1(\tau_f)-\A^{\tau_f}\1_0$.

\subsection{Output controllable system}\label{outcon}
Consider the output controllable system (\ref{sys2}) with $\1_0=\1(0)$ and $\y_f=\y(\tau_f)$.
Let $\C=[\mathbf{I}_r, \mathbf{0}]\in \R^{r\times n}$ with $\mathbf{I}_r$ being $r$-order identity matrix and rank~$\mathcal{C}(\A, \B, \C)=r$ holds.
Analogously, construct Hamilton function as Eq.~(\ref{hm}).
Then Eqs.~(\ref{ze2})\,(\ref{fopu}) and (\ref{opu1}) hold.
From the iterative equation (\ref{ze2}), we have
$
(\A^{\text{T}})^\tau\lambda(\tau)=\lambda(0),
$
which further leads to
$
\lambda(0)=(\A^{\text{T}})^{\tau_f}\lambda(\tau_f).
$
Moreover, let $\lambda(\tau_f)=\C^{\text{T}}\hat{\lambda}_f$, and we have
\begin{equation}\label{Aeq9}
\lambda(\tau)=(\A^{\text{T}})^{\tau_f-\tau}\lambda(\tau_f)=(\A^{\text{T}})^{\tau_f-\tau}\C^{\text{T}}\hat{\lambda}_f.
\end{equation}
And then, substituting Eqs.~(\ref{fopu}) and (\ref{opu1}) into (\ref{Aeq2}), we have
\begin{equation}\notag
\1(\tau)=\A^\tau\1_0-\sum^{\tau-1}_{i=0}\A^{\tau-i-1}\B\B^{\text{T}}\lambda(i+1)
=\A^\tau\1_0-\sum^{\tau-1}_{i=0}\A^{\tau-i-1}\B\B^{\text{T}}(\A^{\text{T}})^{\tau-i-1}\C^{\text{T}}\hat{\lambda}_f.
\end{equation}
Since the controllability Gramian matrix of system (\ref{sys1}) is $\W$,
we  obtain
\begin{equation*}
\begin{cases}
\1(\tau_f)=\A^{\tau_f}\1_0-\W\C^{\text{T}}\hat{\lambda}_f,\\
\y(\tau_f)=\C\A^{\tau_f}\1_0-\C\W\C^{\text{T}}\hat{\lambda}_f,
\end{cases}
\end{equation*}
and the second equality is equivalent to
\begin{equation}\label{Aeq10}
\hat{\lambda}_f=-(\C\W\C^\text{T})^{-1}(\y_f-\C\A^{\tau_f}\1_0).
\end{equation}
Then the optimal input (\ref{opu1}) with (\ref{Aeq9}) and (\ref{Aeq10}) is
\begin{equation*}\label{opu2}
\U^*(\tau)=\B^{\text{T}}(\A^{\text{T}})^{\tau_f-\tau-1}\C^{\text{T}}(\C\W\C^\text{T})^{-1}(\y_f-\C\A^{\tau_f}\1_0).
\end{equation*}
Denoting $\beta=\y_f-\C\A^{\tau_f}\1_0$, the minimum energy cost is
\begin{align*}
E&=\sum^{\tau_f-1}_{k=0}\U^{\text{T}}(k)\U(k)\notag\\
&=\beta^{\text{T}}[(\C\W\C^{\text{T}})^{-1}]^{\text{T}}\C\underbrace{\sum^{\tau_f-1}_{k=0}\A^{\tau_f-k-1}\B\B^{\text{T}} (\A^{\text{T}})^{\tau_f-k-1}}_{\W}\C^{\text{T}} (\C\W\C^{\text{T}})^{-1}\beta\notag\\
&=(\y_f-\C\A^{\tau_f}\1_0)^{\text{T}}(\C\W\C^{\text{T}})^{-1}(\y_f-\C\A^{\tau_f}\1_0).\label{minEout}
\end{align*}

\section{Energy scaling}
\subsection{Energy scaling for full controllability}\label{Energycom}
In this subsection, we assume system (\ref{sys1}) is fully controllable, which implies the Gramian matrix $\W$ is invertible.
Next, we give detailed analysis on approximations of the minimum and maximum eigenvalues of $\W$.
%\begin{equation}\label{upM}
%\lambda_{\max}(\W)\approx f(\overline{\alpha}, \overline{\beta})
%\end{equation}
%\begin{equation}\label{lowM}
%\lambda_{\min}(\W)\approx \frac{1}{f(\underline{\alpha}, \underline{\beta})}.
%\end{equation}
For $\A=\A^{\T}$, we assume that its eigenvalues satisfy $|\lambda_1|\leq|\lambda_2|\leq\cdots\leq|\lambda_n|$.
Via the orthogonal decomposition $\A=\p\la\p^{\T}$ with $\p=(p_{ij})_{nn}$, then Gramian matrix $\W$ is equivalent to
\begin{equation*}\label{WT}
%\begin{split}
\W=\sum^{\tau_f-1}_{\tau=0}(\p\la\p^{\T})^\tau\B\B^{\T}(\p\la\p^{\T})^\tau
%&=\sum^{\tau_f-1}_{\tau=0}\p\la^\tau\p^{\T}\B\B^{\T}\p\la^\tau\p^{\T}\\
=\p\sum^{\tau_f-1}_{\tau=0}\la^\tau\p^{\T}\B\B^{\T}\la^\tau\p^{\T}.
%\end{split}.
\end{equation*}
By letting
%\begin{equation*}\label{M}
$\M=\sum^{\tau_f-1}_{\tau=0}\la^\tau\p^{\T}\B\B^{\T}\p\la^\tau=\sum^{\tau_f-1}_{\tau=0}\M_\tau,$
%\end{equation*}
we have $\W=\p\M\p^{\T}$, which indicates all eigenvalues of $\M$ and $\W$ are the same.
Furthermore, $\M_{\tau}$ has the following form
\begin{equation}\notag
\begin{split}
\M_{\tau}&=
\begin{bmatrix}
\lambda^{\tau}_1&&&\\
&\lambda^{\tau}_2&&\\
&&\ddots&\\
&&&\lambda^{\tau}_n
\end{bmatrix}
\begin{bmatrix}
q_{11}&q_{12}&\cdots&q_{1n}\\
q_{21}&q_{22}&\cdots&q_{2n}\\
\vdots&\vdots&\ddots&\vdots\\
q_{n1}&q_{n2}&\cdots&q_{nn}
\end{bmatrix}
\begin{bmatrix}
\lambda^{\tau}_1&&&\\
&\lambda^{\tau}_2&&\\
&&\ddots&\\
&&&\lambda^{\tau}_n
\end{bmatrix}\\
&\mathop{=}\limits^{(i,j)}q_{ij}\lambda^{\tau}_i\lambda^{\tau}_j,
\end{split}
\end{equation}
where $\Q=\p^{\T}\B\B^{\T}\p=(q_{ij})_{nn}$.
Thus, we have
\begin{equation}\label{Mij}
\M(i, j)=\sum^{\tau_f-1}_{t=0}q_{ij}\lambda^t_i\lambda^t_j=q_{ij}\sum^{\tau_f-1}_{t=0}\lambda^t_i\lambda^t_j=q_{ij} \frac{1-(\lambda_i\lambda_j)^{\tau_f}}{1-\lambda_i\lambda_j}.
\end{equation}
Note that, we have
$
\M(i, j)=q_{ij}\tau_f~\text{if}~ \lambda_i\lambda_j=1,$
since
$
\lim\limits_{\lambda_i\lambda_j\rightarrow1}\frac{1-(\lambda_i\lambda_j)^{\tau_f}}{1-\lambda_i\lambda_j}=\frac{-\tau_f(\lambda_i\lambda_j)^{\tau_f-1}}{-1}=\tau_f.
$
Therefore, in following analysis, we consider the form of $\M$ as Eq.~(\ref{Mij}).

\subsubsection{$n$ driver nodes}
In the case of $n$ driver nodes, each node  receives an independent input signal separately.
Then the corresponding matrix $\Q=\I_n$ causes $\M$ to be a diagonal matrix
$$
\M=
\begin{bmatrix}
\frac{1-\lambda^{2\tau_f}_1}{1-\lambda^2_1}&&\\
&\ddots&\\
&&\frac{1-\lambda^{2\tau_f}_n}{1-\lambda^2_n}
\end{bmatrix}.
$$
Apparently, the function
$$
y=f(x)=\frac{1-x^{2\tau_f}}{1-x^2}
$$
is the sum of geometric sequences
$$
y=f(x)=\sum^{\tau_f-1}_{t=0}x^t=1+x^2+x^4+x^6+\cdots+x^{2(\tau_f-1)},
$$
which is an even function with $\dot{f}(x)>0 ~\text{for}~ x>0$, and $ \dot{f}(x)=0$ at $x=0$.
Therefore, the function $f(x)$  increases as the variable $|x|$ increases.
And then, the minimum and maximum eigenvalues of $\M$ is $\frac{1-\lambda^{2\tau_f}_1}{1-\lambda^2_1}$ and $\frac{1-\lambda^{2\tau_f}_n}{1-\lambda^2_n}$, respectively.
In addition, when time scale $\tau_f$ is large, for $\frac{1-\lambda^{2\tau_f}_1}{1-\lambda^2_1}$ and $\frac{1-\lambda^{2\tau_f}_n}{1-\lambda^2_n}$, we have
\begin{equation*}
\lambda_{\max}(\M)
\begin{cases}
\approx \frac{1}{1-\lambda^2_n},~&\text{if}~ |\lambda_n|<1;\\
\approx \tau_f,~&\text{if}~ |\lambda_n|=1;\\
\thicksim \lambda_n^{2\tau_f-2}, ~&\text{if}~ |\lambda_n|>1,
\end{cases}
\end{equation*}
and
\begin{equation*}
\lambda_{\min}(\M)
\begin{cases}
\approx \frac{1}{1-\lambda^2_1},~&\text{if}~ |\lambda_1|<1;\\
\approx \tau_f,~&\text{if}~ |\lambda_1|=1;\\
\thicksim \lambda_1^{2\tau_f-2}, ~&\text{if}~ |\lambda_1|>1.
\end{cases}
\end{equation*}

\subsubsection{$1$ driver node}\label{1drivecom}
In the case of $1$ driver node, we assume the sole node $h$ as the driver node.
Then, we have the specific form of $\M$ as $\M(i, j)=p_{hi}p_{hj}\frac{1-(\lambda_i\lambda_j)^{\tau_f}}{1-\lambda_i\lambda_j}$. % with $q_{ij}=p_{hi}p_{hj}$.
For the specific form of $\M^2$, we get
$$
(\M^2)_{i, j}=\sum^n_{k=1}p_{hi}p_{hj}p^2_{hk}\frac{1-(\lambda_i\lambda_k)^{\tau_f}}{1-\lambda_i\lambda_k}~\frac{1-(\lambda_k\lambda_j)^{\tau_f}}{1-\lambda_k\lambda_j}.
$$
According to the definitions of $\overline{\alpha}$ and $\overline{\beta}$, we have
\begin{equation*}\label{alphaov}
\overline{\alpha}=\sum_{i=1}^n\sum^n_{k=1}p_{hi}^2p^2_{hk}\left(\frac{1-(\lambda_i\lambda_k)^{\tau_f}}{1-\lambda_i\lambda_k}\right)^2
\end{equation*}
and
\begin{equation*}\label{betaov}
\overline{\beta}=\sum^n_{j=1}\sum^n_{i=1}\left(\sum^n_{k=1}p_{hi}p_{hj}p^2_{hk}\frac{1-(\lambda_i\lambda_k)^{\tau_f}}{1- \lambda_i\lambda_k}~\frac{1-(\lambda_k\lambda_j)^{\tau_f}}{1-\lambda_k\lambda_j}\right)^2.
\end{equation*}
More specifically, when the maximum eigenvalue of $\A$ satisfies $|\lambda_n|< 1$, with the approximation of $1-(\lambda_i\lambda_j)^{\tau_f}\approx 1$ for $i, j=1, 2, \dots, n$, we have
\begin{equation}\label{alphaup1}
\overline{\alpha}\approx\sum_{i=1}^n\sum^n_{k=1}p_{hi}^2p^2_{hk}\left(\frac{1}{1-\lambda_i\lambda_k}\right)^2
\end{equation}
and
\begin{equation}\label{betaup1}
\overline{\beta}\approx\sum^n_{j=1}\sum^n_{i=1}\left(\sum^n_{k=1}p_{hi}p_{hj}p^2_{hk}\frac{1}{1-\lambda_i\lambda_k}~\frac{1}{1-\lambda_k\lambda_j}\right)^2.
\end{equation}
When the maximum eigenvalue of $\A$ satisfies $|\lambda_n| = 1$, with the approximation of $\frac{1-(\lambda_i\lambda_j)^{\tau_f}}{1-(\lambda_i\lambda_j)}\approx \tau_f$ for $i, j=1, 2, \dots, n$, we have
\begin{equation*}\label{alphaup2}
\overline{\alpha}=\sum_{i=1}^n\sum^n_{k=1}p_{hi}^2p^2_{hk}\left(\frac{1-(\lambda_i\lambda_k)^{\tau_f}}{1-\lambda_i\lambda_k}\right)^2\thicksim {\tau_f}^2
\end{equation*}
and
\begin{equation*}\label{betaup2}
\overline{\beta}=\sum^n_{j=1}\sum^n_{i=1}\left(\sum^n_{k=1}p_{hi}p_{hj}p^2_{hk}\frac{1-(\lambda_i\lambda_k)^{\tau_f}}{1- \lambda_i\lambda_k}~\frac{1-(\lambda_k\lambda_j)^{\tau_f}}{1-\lambda_k\lambda_j}\right)^2\thicksim \tau_f^4.
\end{equation*}
When the maximum eigenvalue of $\A$ satisfies $|\lambda_n| > 1$, we have
\begin{equation*}\label{alphaup3}
\overline{\alpha}=\sum_{i=1}^n\sum^n_{k=1}p_{hi}^2p^2_{hk}\left(\frac{1-(\lambda_i\lambda_k)^{\tau_f}}{1-\lambda_i\lambda_k}\right)^2\thicksim \lambda_n^{4\tau_f-4}
\end{equation*}
and
\begin{equation*}\label{betaup3}
\overline{\beta}=\sum^n_{j=1}\sum^n_{i=1}\left(\sum^n_{k=1}p_{hi}p_{hj}p^2_{hk}\frac{1-(\lambda_i\lambda_k)^{\tau_f}}{1- \lambda_i\lambda_k}~\frac{1-(\lambda_k\lambda_j)^{\tau_f}}{1-\lambda_k\lambda_j}\right)^2\thicksim \lambda_n^{8\tau_f-8}.
\end{equation*}

In order to get values of $\underline{\alpha}$ and $\underline{\beta}$, the key precondition relies on $\M^{-1}$.
When all eigenvalues of $\A$ satisfy $|\lambda_i|<1$, we have $\M(i, j)=\frac{p_{hi}p_{hj}}{1-\lambda_i\lambda_j}$ with the approximation of $1-(\lambda_i\lambda_j)^{\tau_f}\approx 1$.
Furthermore, the corresponding inverse matrix $\M^{-1}$ is
\begin{equation}\notag
\M^{-1}(i, j)=
\begin{cases}
\frac{(1-\lambda_i^2)\prod\limits_{k\neq i}(1-\lambda_i\lambda_k)^2}{\prod\limits_{k\neq i}(\lambda_i-\lambda_k)^2p^2_{hi}}, &i=j;\\
\frac{(1-\lambda_i^2)(1-\lambda_j^2)\prod\limits_{k< l}\prod\limits_{l=2}^{n}(1-\lambda_k\lambda_l)}{(\lambda_i-\lambda_j)^2\prod\limits_{k\neq i,j}(\lambda_i-\lambda_k)(\lambda_j-\lambda_k)p_{hi}p_{hj}}, &i\neq j.
\end{cases}
\end{equation}
Hence, in this case, matrix $\M^2$ has the following form
\begin{equation}\notag
\begin{split}
(\M^2)_{i, j}=&\sum^n \limits_{k=1, k\neq i, j}\M(i, k)\M(k, j)+\M(i, i)\M(i, j)+\M(i, j)\M(j, j)\\
=&\sum^n_{k=1, k\neq i, j}\frac{(1-\lambda_i^2)(1-\lambda_j^2)(1-\lambda_k^2)^2\prod\limits_{r< l}\prod\limits_{l=2}^{n}(1-\lambda_r\lambda_l)^2}{(\lambda_i-\lambda_k)^2(\lambda_j-\lambda_k)^2\prod\limits_{r\neq i,k}\prod \limits_{l\neq k, j}(\lambda_i-\lambda_r)(\lambda_k-\lambda_r)(\lambda_k-\lambda_l)(\lambda_j-\lambda_l)p_{hi}p_{hj}p^2_{hk}}\\
&+\frac{(1-\lambda_i^2)\prod\limits_{r\neq i}(1-\lambda_i\lambda_r)^2}{\prod\limits_{r\neq i}(\lambda_i-\lambda_r)^2p^2_{hi}} \frac{(1-\lambda_i^2)(1-\lambda_j^2)\prod\limits_{r< l}\prod\limits_{l=2}^{n}(1-\lambda_r\lambda_l)}{(\lambda_i-\lambda_j)^2\prod\limits_{r\neq i,j}(\lambda_i-\lambda_r)(\lambda_j-\lambda_r)p_{hi}p_{hj}}\\
&+\frac{(1-\lambda_j^2)\prod\limits_{r\neq j}(1-\lambda_j\lambda_r)^2}{\prod\limits_{r\neq j}(\lambda_j-\lambda_r)^2p^2_{hj}} \frac{(1-\lambda_i^2)(1-\lambda_j^2)\prod\limits_{r< l}\prod\limits_{l=2}^{n}(1-\lambda_r\lambda_l)}{(\lambda_i-\lambda_j)^2\prod\limits_{r\neq i,j}(\lambda_i-\lambda_r)(\lambda_j-\lambda_r)p_{hi}p_{hj}}.
\end{split}
\end{equation}
From that, we can derive $\underline{\alpha}$ and $\underline{\beta}$ as
\begin{equation*}\label{alphalow1}
\underline{\alpha}=\sum^n_{i=1}\left(\frac{(1-\lambda_i^2)\prod\limits_{k\neq i}(1-\lambda_i\lambda_k)^2}{\prod\limits_{k\neq i}(\lambda_i-\lambda_k)^2p^2_{hi}}\right)^2+
\sum^n_{j=1}\sum^n_{i=1, i\neq j}\left(\frac{(1-\lambda_i^2)(1-\lambda_j^2)\prod\limits_{k< l}\prod\limits_{l=2}^{n}(1-\lambda_k\lambda_l)}{(\lambda_i-\lambda_j)^2\prod\limits_{k\neq i,j}(\lambda_i-\lambda_k)(\lambda_j-\lambda_k)p_{hi}p_{hj}}\right)^2
\end{equation*}
and
\begin{align*}
\underline{\beta}=&\sum^n_{j=1}\sum^n_{i=1}\left(\sum^n_{k=1, k\neq i, j}\frac{(1-\lambda_i^2)(1-\lambda_j^2)(1-\lambda_k^2)^2\prod\limits_{r< l}\prod\limits_{l=2}^{n}(1-\lambda_r\lambda_l)^2}{(\lambda_i-\lambda_k)^2(\lambda_j-\lambda_k)^2\prod \limits_{r\neq i,k; l\neq k, j}(\lambda_i-\lambda_r)(\lambda_k-\lambda_r)(\lambda_k-\lambda_l)(\lambda_j-\lambda_l)p_{hi}p_{hj}p^2_{hk}}\right.\notag\\
&\left.+\frac{(1-\lambda_i^2)\prod\limits_{r\neq i}(1-\lambda_i\lambda_r)^2}{\prod\limits_{r\neq i}(\lambda_i-\lambda_r)^2p^2_{hi}} \frac{(1-\lambda_i^2)(1-\lambda_j^2)\prod\limits_{r< l}\prod\limits_{l=2}^{n}(1-\lambda_r\lambda_l)}{(\lambda_i-\lambda_j)^2\prod\limits_{r\neq i,j}(\lambda_i-\lambda_r)(\lambda_j-\lambda_r)p_{hi}p_{hj}}\right.\notag\\
&\left.+\frac{(1-\lambda_j^2)\prod\limits_{r\neq j}(1-\lambda_j\lambda_r)^2}{\prod\limits_{r\neq j}(\lambda_j-\lambda_r)^2p^2_{hj}} \frac{(1-\lambda_i^2)(1-\lambda_j^2)\prod\limits_{r< l}\prod\limits_{l=2}^{n}(1-\lambda_r\lambda_l)}{(\lambda_i-\lambda_j)^2\prod\limits_{r\neq i,j}(\lambda_i-\lambda_r)(\lambda_j-\lambda_r)p_{hi}p_{hj}}\right)^2.\label{betalow1}
\end{align*}
When $|\lambda_1|<1$ and other eigenvalues satisfy $|\lambda_1|\leq\cdots\leq|\lambda_l|<1$, $|\lambda_{l+1}|=\cdots=|\lambda_{l+r}|=1$, $1<|\lambda_{l+r+1}|\leq\cdots\leq|\lambda_n|$, the corresponding $\M$ is
\begin{equation*}\label{MTP}
\M=
\begin{bmatrix}
\M_{1}&\M_{2}&\M_{3}\\
\M^{\T}_{2}&\M_{4}&\M_{5}\\
\M^{\T}_{3}&\M^{\T}_{5}&\M_{6}
\end{bmatrix}
\end{equation*}
where
$$
\M_{1}=
\begin{bmatrix}
\frac{p^2_{h1}}{1-\lambda_1^2}&\cdots&\frac{p_{h1}p_{hl}}{1-\lambda_1\lambda_l}\\
\vdots&\ddots&\vdots\\
\frac{p_{h1}p_{hl}}{1-\lambda_1\lambda_l}&\cdots&\frac{p^2_{hl}}{1-\lambda_l^2}
\end{bmatrix},\quad
\M_{2}=
\begin{bmatrix}
\frac{p_{h1}p_{l+1} (1-(\lambda_1\lambda_{l+1})^{\tau_f})}{1-\lambda_1\lambda_{l+1}}&\cdots&\frac{p_{h1}p_{l+r} (1-(\lambda_1\lambda_{l+r})^{\tau_f})}{1-\lambda_1\lambda_{l+r}}\\
\vdots&\ddots&\vdots\\
\frac{p_{hl}p_{l+1} (1-(\lambda_l\lambda_{l+1})^{\tau_f})}{1-\lambda_l\lambda_{l+1}}&\cdots&\frac{p_{hl}p_{l+r} (1-(\lambda_l\lambda_{l+r})^{\tau_f})}{1-\lambda_l\lambda_{l+r}}
\end{bmatrix},
$$
$$
\M_{3}=
\begin{bmatrix}
\frac{p_{h1}p_{l+r+1} (1-(\lambda_1\lambda_{l+r+1})^{\tau_f})}{1-\lambda_1\lambda_{l+r+1}}&\cdots&\frac{p_{h1}p_{n} (1-(\lambda_1\lambda_{n})^{\tau_f})}{1-\lambda_1\lambda_{n}}\\
\vdots&\ddots&\vdots\\
\frac{p_{hl}p_{l+r+1} (1-(\lambda_l\lambda_{l+r+1})^{\tau_f})}{1-\lambda_l\lambda_{l+r+1}}&\cdots&\frac{p_{hl}p_{n} (1-(\lambda_l\lambda_{n})^{\tau_f})}{1-\lambda_l\lambda_{n}}
\end{bmatrix},$$$$
\M_{4}=
\begin{bmatrix}
p_{h\,l+1}^2{\tau_f}&\cdots&p_{h\,l+1}p_{h\,l+r}{\tau_f}\\
\vdots&\ddots&\vdots\\
p_{h\,l+1}p_{h\,l+r}{\tau_f}&\cdots &p_{h\,l+r}^2{\tau_f}
\end{bmatrix},
$$
$$
\M_{5}=
\begin{bmatrix}
\frac{p_{h\,l+1}p_{h\,l+r+1} (1-(\lambda_{h\, l+1}\lambda_{l+r+1})^{\tau_f})}{1-\lambda_{l+1}\lambda_{l+r+1}}&\cdots&\frac{p_{h\, l+1}p_{h\,n} (1-(\lambda_{l+1}\lambda_{n})^{\tau_f})}{1-\lambda_{l+1}\lambda_{n}}\\
\vdots&\ddots&\vdots\\
\frac{p_{h\,l+r}p_{l+r+1} (1-(\lambda_{l+r}\lambda_{l+r+1})^{\tau_f})}{1-\lambda_{l+r}\lambda_{l+r+1}}&\cdots&\frac{p_{h\,l+r}p_{h\,n} (1-(\lambda_{l+r}\lambda_{n})^{\tau_f})}{1-\lambda_{l+r}\lambda_{n}}
\end{bmatrix},
$$
$$
\M_{6}=
\begin{bmatrix}
\frac{p_{h\,l+r+1}^2 (1-\lambda_{l+r+1}^{2{\tau_f}})}{1-\lambda_{l+r+1}^2}&\cdots&\frac{p_{h\, l+r+1}p_{h\,n} (1-(\lambda_{l+r+1}\lambda_{n})^{\tau_f})}{1-\lambda_{l+r+1}\lambda_{n}}\\
\vdots&\ddots&\vdots\\
\frac{p_{h\, l+r+1}p_{h\,n} (1-(\lambda_{l+r+1}\lambda_{n})^{\tau_f})}{1-\lambda_{l+r+1}\lambda_{n}}&\cdots&\frac{p_{hn}^2 (1-\lambda_{n}^{2{\tau_f}})}{1-\lambda_{n}^2}
\end{bmatrix}.
$$
In this case, it is difficult to calculate $\M^{-1}$ directly.
For that, we utilize $\M^{-1}=\frac{\M^*}{|\M|}$ with $\M^*$ being adjoint matrix of $\M$.
Intuitively, it is easy to get
$
|\M|\thicksim {\tau_f}^r~\prod^n_{i={l+r+1}}\lambda_i^{2{\tau_f}}.
$
For adjoint matrix $\M^*$, we have
$$
\M^*(i, j)\thicksim
\begin{cases}
{\tau_f}^r\cdot \prod^n_{i={l+r+1}}\lambda_i^{2{\tau_f}}, &i, j\leq l;\\
{\tau_f}^a\cdot \prod^n_{i={l+r+1}}\lambda_i^{2{\tau_f}}/k, &\text{otherwise},
\end{cases}
$$
where $a=r-1 ~\text{or}~ k\neq 1 ~\text{with}~ k=(\lambda_{l_1}\lambda_{l_2})^{{\tau_f}}$, $l_1, l_2\geq l+r+1$.
Therefore, we get $\M^{-1}$ as
$$
\M^{-1}(i,j)\approx
\begin{cases}
c_{ij}\neq 0, &i, j\leq l;\\
0, &\text{otherwise}.
\end{cases}
$$
Then in the process of calculating $\underline{\alpha}$ and $\underline{\beta}$, we find that elements of the first $l$ rows and $l$ columns of $\M^{-1}$ dominate, i.e., $c_{ij}, i, j\leq l$ are adequate.
Thus, in order to get the specific forms of $\underline{\alpha}$ and $\underline{\beta}$, we employ the inverse matrix of $\M_{1}$ to replace the inverse matrix of $\M$.
And then, the corresponding $\underline{\alpha}$ and $\underline{\beta}$ are
\begin{equation}\label{alphalow1}
\underline{\alpha}=\sum^l_{i=1}\left(\frac{(1-\lambda_i^2)\prod\limits_{k\neq i}(1-\lambda_i\lambda_k)^2}{\prod\limits_{k\neq i}(\lambda_i-\lambda_k)^2p^2_{hi}}\right)^2+
\sum^l_{j=1}\sum^l_{i=1, i\neq j}\left(\frac{(1-\lambda_i^2)(1-\lambda_j^2)\prod\limits_{k< d}\prod\limits_{d=2}^{n}(1-\lambda_k\lambda_d)}{(\lambda_i-\lambda_j)^2\prod\limits_{k\neq i,j}(\lambda_i-\lambda_k)(\lambda_j-\lambda_k)p_{hi}p_{hj}}\right)^2
\end{equation}
and
\begin{equation}\label{betalow1}
\begin{split}
\underline{\beta}=&\sum^l_{j=1}\sum^l_{i=1}\left(\sum^l_{k=1, k\neq i, j}\frac{(1-\lambda_i^2)(1-\lambda_j^2)(1-\lambda_k^2)^2\prod\limits_{r< d}\prod\limits_{d=2}^{l}(1-\lambda_r\lambda_d)^2}{(\lambda_i-\lambda_k)^2(\lambda_j-\lambda_k)^2\prod \limits_{r\neq i,k; d\neq k, j}(\lambda_i-\lambda_r)(\lambda_k-\lambda_r)(\lambda_k-\lambda_d)(\lambda_j-\lambda_d)p_{hi}p_{hj}p^2_{hk}}\right.\\
&+\frac{(1-\lambda_i^2)\prod\limits_{r\neq i}(1-\lambda_i\lambda_r)^2}{\prod\limits_{r\neq i}(\lambda_i-\lambda_r)^2p^2_{hi}} \frac{(1-\lambda_i^2)(1-\lambda_j^2)\prod\limits_{r< d}\prod\limits_{d=2}^{n}(1-\lambda_r\lambda_d)}{(\lambda_i-\lambda_j)^2\prod\limits_{r\neq i,j}(\lambda_i-\lambda_r)(\lambda_j-\lambda_r)p_{hi}p_{hj}}\\
&\left.+\frac{(1-\lambda_j^2)\prod\limits_{r\neq j}(1-\lambda_j\lambda_r)^2}{\prod\limits_{r\neq j}(\lambda_j-\lambda_r)^2p^2_{hj}} \frac{(1-\lambda_i^2)(1-\lambda_j^2)\prod\limits_{r< d}\prod\limits_{d=2}^{n}(1-\lambda_r\lambda_d)}{(\lambda_i-\lambda_j)^2\prod\limits_{r\neq i,j}(\lambda_i-\lambda_r)(\lambda_j-\lambda_r)p_{hi}p_{hj}}\right)^2.
\end{split}
\end{equation}

When $|\lambda_1|=1$, i.e., $
\M=
\begin{bmatrix}
\M_{4}&\M_{5}\\
\M^{\T}_{5}&\M_{6}
\end{bmatrix}
$
with $l=0$ in (\ref{MTP}), for $\M^{-1}=\frac{\M^*}{|\M|}$, we have
$
|\M|\thicksim {\tau_f}^r\cdot \prod^n_{i={r+1}}\lambda_i^{2{\tau_f}}
$
and
$$
\M^*(i, j)\thicksim
\begin{cases}
{\tau_f}^{r-1}\cdot \prod^n_{i={r+1}}\lambda_i^{2{\tau_f}}, &i, j\leq r;\\
{\tau_f}^{a}\cdot \prod^n_{i={r+1}}\lambda_i^{2{\tau_f}}/k, &\text{otherwise},
\end{cases}
$$
$\text{where}~ a=r-1  ~\text{or}~a=r~\&~ k=(\lambda_{l_1}\lambda_{l_2})^{{\tau_f}},$
which lead to
$$
\M^{-1}(i,j)\thicksim
\begin{cases}
{\tau_f}^{-1}, &i, j\leq r;\\
{\tau_f}^{b}(\lambda_{l_1}\lambda_{l_2})^{-{\tau_f}},&\text{otherwise} (\text{with}~b=0~\text{or}-1).
\end{cases}
$$
Due to that $
\lim_{{\tau_f}\rightarrow \infty}\frac{{\tau_f}^{b}(\lambda_{l_1}\lambda_{l_2})^{-{\tau_f}}}{{\tau_f}^{-1}}=0
$, we have
$$\underline{\alpha}\thicksim {\tau_f}^{-2}\quad \text{and} \quad \underline{\beta}\thicksim {\tau_f}^{-4}.$$

When $|\lambda_1|>1$,  i.e., $
\M=\M_{6}
$
with $l=r=0$ in (\ref{MTP}), for $\M^{-1}=\frac{\M^*}{|\M|}$, we have
$
|\M|\thicksim \prod^n_{i={1}}\lambda_i^{2{\tau_f}}=|\A|^{2{\tau_f}}
$
and
$
\M^*(i,j)\thicksim\frac{\prod^n_{i={1}}\lambda_i^{2{\tau_f}}}{\lambda_i^{\tau_f}\lambda^{\tau_f}_j},
$
which lead to
$\M^{-1}(i, j)\thicksim(\lambda_i\lambda_j)^{-{\tau_f}}.$
Moreover, in calculating  $\underline{\alpha}$ and $\underline{\beta}$, $\lambda_1^{-2{\tau_f}}$ dominates.
Therefore, we have
$$\underline{\alpha}\thicksim \lambda_1^{-4{\tau_f}},\quad\text{and} \quad \underline{\beta}\thicksim \lambda_1^{-8{\tau_f}}.$$

\subsubsection{$m$ driver nodes}
In the case of $m$ driver nodes, the indexes of driver nodes can be denoted by $d_1, d_2, \cdots, d_m$.
Then the corresponding input matrix is $\B=[e_{d_1}, e_{d_2}, \cdots, e_{d_m}]\in \R^{n\times m}$ with $e_i$ being the $i$th column of identity matrix.
Accordingly, we have $\M({i, j})=q_{ij}\frac{1-(\lambda_i\lambda_j)^{\tau_f}}{1-\lambda_i\lambda_j}$ with $q_{ij}=\sum^m_{k=1}p_{d_k i}p_{d_k j}$.
In the analysis of $1$ driver \text{node}, we find that the form of $q_{ij}$ has no essential effect on the main analysis process.
For example, when $|\lambda_n|<1$, we have $\M({i, j})\approx \sum^m_{k=1}p_{d_k i}p_{d_k j}\frac{1}{1-\lambda_i\lambda_j}$ and
$\M^2(i, j)=\sum_{l=1}^n(\sum^m_{k=1}p_{d_k i}p_{d_k l}\frac{1}{1-\lambda_i\lambda_l})(\sum^m_{k=1}p_{d_k l}p_{d_k j}\frac{1}{1-\lambda_l\lambda_j})$.
Furthermore, we have the corresponding $\overline{\alpha}$ and $\overline{\beta}$
\begin{equation}\label{alphaup4}
\overline{\alpha}=\sum_{j=1}^n\sum_{i=1}^n\left(\frac{\sum^m_{k=1}p_{d_k i}p_{d_k j}}{1-\lambda_i\lambda_j}\right)^2
\end{equation}
and
\begin{equation}\label{betaup4}
\overline{\beta}=\sum_{j=1}^n\sum_{i=1}^n\left(
\sum_{l=1}^n\left(\sum^m_{k=1}p_{d_k i}p_{d_k l}\frac{1}{1-\lambda_i\lambda_l}\right)\left(\sum^m_{k=1}p_{d_k l}p_{d_k j}\frac{1}{1-\lambda_l\lambda_j}\right)   \right)^2.
\end{equation}
Analogously, other cases can also be derived and thus omitted here.

In summary, in different cases of $\A$ with different properties, parameters $\overline{\alpha}, \overline{\beta}, \underline{\alpha}$ and $\underline{\beta}$ are obtained. And the corresponding bounds are acquired accordingly as shown in table~\ref{table2}.
%In addition, the corresponding numerical verification is presented in fig.~\ref{figure5}.

\begin{table}[!http]
\centering\caption{Lower and upper bounds of the minimum energy for a fully controllable network.
}
\fontsize{8}{15}\selectfont
\begin{tabular}{cc|c|c|c}
\toprule[2pt]
\multicolumn{2}{c|}{Number of driver nodes}&$1$&$m (m<n)$&$n$\\
\toprule[1pt]
%\hline
\multirow{3}{*}{Lower bound $\underline{E}$}&$|\lambda_{n}|<1$&Eq.(\ref{lowE}) with (\ref{alphaup1})(\ref{betaup1})&Eq.(\ref{lowE}) with (\ref{alphaup4})(\ref{betaup4})&$1-\lambda_n^2$\\
%\hline
&$|\lambda_{n}|=1$&$ \sim \tau_f^{-1}$& $\sim \tau_f^{-1}$& $\sim \tau_f^{-1}$\\
%\hline
&$|\lambda_{n}|>1$&$\sim \lambda_{n}^{2-2\tau_f}$&$\sim \lambda_{n}^{2-2\tau_f}$&$\sim \lambda_{n}^{2-2\tau_f}$\\
\toprule[1pt]
\multirow{3}{*}{Upper bound $\overline{E}$}&$|\lambda_{1}|<1$&Eq.(\ref{upE}) with (\ref{alphalow1})(\ref{betalow1})&constant&$1-\lambda_1^2$\\
%\hline
&$|\lambda_1|=1$&$ \sim \tau_f^{-1}$& $\sim \tau_f^{-1}$& $\sim \tau_f^{-1}$\\
%\hline
&$|\lambda_1|>1$&$\sim \lambda_{1}^{-2\tau_f}$&$\sim \lambda_{1}^{-2\tau_f}$&$\sim \lambda_{1}^{-2\tau_f}$\\
\toprule[2pt]
\end{tabular}
\label{table2}
\end{table}
%\begin{figure}[htbp]
%\centering
%\begin{minipage}[t]{1\textwidth}
%\centering
%\includegraphics[width=12cm]{../Public/f262.pdf}
%\end{minipage}
%\caption{The lower and upper bounds of energy for full control.
%Here, we adopt scale free networks for numerical verification.
%All paraments selected are the same as those in fig.~\ref{figure7}.
%Therefore, eigenvalues of matrix $\A$ can present different characters.
%Then, we need to select a driver node such that the network is fully controllable.
%Similarly,  in each subplot, it is obvious that the generated pattern of analytical derivation almost overlaps that of numerical calculation.
%}\label{figure5}
%\end{figure}

\subsection{Energy scaling for target control}\label{energyuncom}
The essential procedure is to get the minimum and maximum eigenvalues of $\W_{\text{C}}$ in Eq.~(\ref{WCT}).
Furthermore, we employ $\lambda_{\max}(\W)\approx f(\overline{\alpha}, \overline{\beta})$ and $\lambda_{\min}(\W)\approx \frac{1}{f(\underline{\alpha}, \underline{\beta})}$ to approximate the corresponding eigenvalues.
Note that subsystem (\ref{eq2}) is controllable.
Analogical to the case of full controllability, we perform the following analysis.
For system (\ref{eq2}), we have $\A_{\text{c}}=\p_{\text{c}}\la_{\text{c}}\p_{\text{c}}^{\T}$ and $\B_{\text{c}}=\2_1\B_r$ with $\la=$diag$(\lambda_{{\text{c}}1}, \lambda_{{\text{c}}2}, \dots, \lambda_{{\text{c}}r})$  and $\B_r$ being the first $r$ rows of $\B$.
It is obvious that $\lambda_{{\text{c}}i}\in \{\lambda_1, \lambda_2, \dots, \lambda_n\}$.
Moreover, the corresponding $\mathcal{W}$ is $\p_{\text{c}}\M_{\text{C}}\p_{\text{c}}^{\T}$, where $\M_{\text{C}}=\sum_{{\tau}=0}^{{\tau_f}-1}\la_{\text{c}}^{\tau}\Q_{\text{c}}\la_{\text{c}}^{\tau}$ with $\Q_{\text{c}}=\p_{\text{c}}^{\T}\B_{\text{c}}\B_{\text{c}}^{\T}\p_{\text{c}}$.
Denoting $\p_R=\2_1^{\T}\p_{\text{c}}$, we have $\Q_{\text{c}}=\p_R^{\T}\B_r\B_r^{\T}\p_R$ with $\Q_{\text{c}}=(q_{ij}^C)_{r\times r}$ and $\p_R=(p_{ij}^R)_{r\times r}$.

In the case of $1$ driver node, %i.e., only the driver node injected by an external input signal, assume the index of the  corresponding driver node as $h$.
the elements of $\M_{\text{C}}$ are $\M_{\text{C}}(i, j)=q_{ij}^C\frac{1-(\lambda_{{\text{c}}i}\lambda_{{\text{c}}j})^{{\tau_f}}}{1-\lambda_{{\text{c}}i}\lambda_{{\text{c}}j}}$ with $q_{ij}^C=p_{hi}^Rp_{hj}^R, i, j=1, 2, \dots, r.$
The specific forms of $\W_{\text{C}}$ and $\W_{\text{C}}^2$ are
\begin{equation*}
\W_{\text{C}}(i, j)=\sum_{k=1}^r\sum_{l=1}^r p_{ik}^Rp_{jl}^Rq_{kl}^C\frac{1-(\lambda_{{\text{c}}k}\lambda_{{\text{c}}l})^{\tau_f}}{1-\lambda_{{\text{c}}k}\lambda_{{\text{c}}l}}
\end{equation*}
and
\begin{equation*}
\W_{\text{C}}^2(i, j)=\sum_{s=1}^r\left(\sum_{k=1}^r\sum_{l=1}^r p_{ik}^Rp_{sl}^Rq_{kl}^C\frac{1-(\lambda_{{\text{c}}k}\lambda_{{\text{c}}l})^{\tau_f}}{1-\lambda_{{\text{c}}k}\lambda_{{\text{c}}l}}\right)
\left( \sum_{k=1}^r\sum_{l=1}^r p_{sk}^Rp_{jl}^Rq_{kl}^C\frac{1-(\lambda_{{\text{c}}k}\lambda_{{\text{c}}l})^{\tau_f}}{1-\lambda_{{\text{c}}k}\lambda_{{\text{c}}l}} \right).
\end{equation*}
According to the definitions of $\overline{\alpha}$ and $\overline{\beta}$, we have
\begin{equation*}
\begin{split}
\overline{\alpha}=\sum_{i=1}^r\sum_{j=1}^r\left( \sum_{k=1}^r\sum_{l=1}^r p_{ik}^Rp_{jl}^Rq_{kl}^C\frac{1-(\lambda_{{\text{c}}k}\lambda_{{\text{c}}l})^{\tau_f}}{1-\lambda_{{\text{c}}k}\lambda_{{\text{c}}l}}\right)^2
\end{split}
\end{equation*}
and
\begin{equation*}
\begin{split}
\overline{\beta}
&=\sum_{i=1}^r\sum_{j=1}^r\left( \sum_{s=1}^r\left(\sum_{k=1}^r\sum_{l=1}^r p_{ik}^Rp_{sl}^Rq_{kl}^C\frac{1-(\lambda_{{\text{c}}k}\lambda_{{\text{c}}l})^{\tau_f}}{1-\lambda_{{\text{c}}k}\lambda_{{\text{c}}l}}\right)
\left( \sum_{k=1}^r\sum_{l=1}^r p_{sk}^Rp_{jl}^Rq_{kl}^C\frac{1-(\lambda_{{\text{c}}k}\lambda_{{\text{c}}l})^{\tau_f}}{1-\lambda_{{\text{c}}k}\lambda_{{\text{c}}l}} \right) \right)^2.
\end{split}
\end{equation*}
Similar to the section~\ref{1drivecom}, based on the approximation of $1-(\lambda_{{\text{c}}i}\lambda_{{\text{c}}j})^{\tau_f}\approx 1$ for $|\lambda_{{\text{c}}r}|<1$, we have
\begin{equation}\label{Aeq11}
\begin{split}
\overline{\alpha}=\sum_{i=1}^r\sum_{j=1}^r\left( \sum_{k=1}^r\sum_{l=1}^r p_{ik}^Rp_{jl}^Rq_{kl}^C\frac{1}{1-\lambda_{{\text{c}}k}\lambda_{{\text{c}}l}}\right)^2
\end{split}
\end{equation}
and
\begin{equation}\label{Aeq12}
\begin{split}
\overline{\beta}
&=\sum_{i=1}^r\sum_{j=1}^r\left( \sum_{s=1}^r\left(\sum_{k=1}^r\sum_{l=1}^r p_{ik}^Rp_{sl}^Rq_{kl}^C\frac{1}{1-\lambda_{{\text{c}}k}\lambda_{{\text{c}}l}}\right)
\left( \sum_{k=1}^r\sum_{l=1}^r p_{sk}^Rp_{jl}^Rq_{kl}^C\frac{1}{1-\lambda_{{\text{c}}k}\lambda_{{\text{c}}l}} \right) \right)^2.
\end{split}
\end{equation}
In the case of $|\lambda_{{\text{c}}r}|=1$, by utilizing $\frac{1-(\lambda_{{\text{c}}i}\lambda_{{\text{c}}j})^{\tau_f}}{1-(\lambda_{{\text{c}}i}\lambda_{{\text{c}}j})}\approx {\tau_f}$, we have
\begin{equation*}
\begin{split}
\overline{\alpha}&\approx\sum_{i=1}^r\sum_{j=1}^r\left( \sum_{k=1}^r\sum_{l=1}^r p_{ik}^Rp_{jl}^Rq_{kl}^C{\tau_f}\right)^2\sim {\tau_f}^2
\end{split}
\end{equation*}
and
\begin{equation*}
\begin{split}
\overline{\beta}
&=\sum_{i=1}^r\sum_{j=1}^r\left( \sum_{s=1}^r\left(\sum_{k=1}^r\sum_{l=1}^r p_{ik}^Rp_{sl}^Rq_{kl}^C{\tau_f}\right)
\left( \sum_{k=1}^r\sum_{l=1}^r p_{sk}^Rp_{jl}^Rq_{kl}^C{\tau_f} \right) \right)^2\sim {\tau_f}^4.
\end{split}
\end{equation*}
In the case of $|\lambda_{{\text{c}}r}|>1$, we have
\begin{equation*}
\begin{split}
\overline{\alpha}=\sum_{i=1}^r\sum_{j=1}^r\left( \sum_{k=1}^r\sum_{l=1}^r p_{ik}^Rp_{jl}^Rq_{kl}^C\frac{1-(\lambda_{{\text{c}}k}\lambda_{{\text{c}}l})^{\tau_f}}{1-\lambda_{{\text{c}}k}\lambda_{{\text{c}}l}}\right)^2\sim \lambda_{{\text{c}}r}^{4{\tau_f}}
\end{split}
\end{equation*}
and
\begin{equation*}
\begin{split}
\overline{\beta}
&=\sum_{i=1}^r\sum_{j=1}^r\left( \sum_{s=1}^r\left(\sum_{k=1}^r\sum_{l=1}^r p_{ik}^Rp_{sl}^Rq_{kl}^C\frac{1-(\lambda_{{\text{c}}k}\lambda_{{\text{c}}l})^{\tau_f}}{1-\lambda_{{\text{c}}k}\lambda_{{\text{c}}l}}\right)
\left( \sum_{k=1}^r\sum_{l=1}^r p_{sk}^Rp_{jl}^Rq_{kl}^C\frac{1-(\lambda_{{\text{c}}k}\lambda_{{\text{c}}l})^{\tau_f}}{1-\lambda_{{\text{c}}k}\lambda_{{\text{c}}l}} \right) \right)^2\\
&\sim \lambda_{{\text{c}}r}^{8{\tau_f}}.
\end{split}
\end{equation*}

To calculate $\underline{\alpha}$ and $\underline{\beta}$, the pivotal is to get $\W_{\text{C}}^{-1}$.
It is clear that $\W_{\text{C}}^{-1}=\2_1^{-1}\mathcal{W}^{-1}(\2_{1}^{\T})^{-1}$ with $\mathcal{W}^{-1}=\p_{\text{c}}\M_{\text{C}}^{-1}\p_{\text{c}}^{\T}$.
Moreover, when $|\lambda_{{\text{c}}i}|<1$, $\M_{\text{C}}(i, j)\approx \frac{p_{hi}^Rp_{hj}^R}{1-\lambda_{{\text{c}}i}\lambda_{{\text{c}}j}}$.
The corresponding elements of  $\M_{\text{C}}^{-1}$ are
\begin{equation}\notag
\M_{\text{C}}^{-1}(i, j)=
\begin{cases}
\frac{(1-\lambda_{{\text{c}}i}^2)\prod\limits_{k\neq i}(1-\lambda_{{\text{c}}i}\lambda_{{\text{c}}k})^2}{\prod\limits_{k\neq i}(\lambda_{{\text{c}}i}-\lambda_{{\text{c}}k})^2(p^{R}_{hi})^2}, &i=j;\\
\frac{(1-\lambda_{{\text{c}}i}^2)(1-\lambda_{{\text{c}}j}^2)\prod\limits_{k< l}\prod\limits_{l=2}^{n}(1-\lambda_{{\text{c}}k}\lambda_{{\text{c}}l})}{(\lambda_{{\text{c}}i}-\lambda_{{\text{c}}j})^2\prod\limits_{k\neq i,j}(\lambda_{{\text{c}}i}-\lambda_{{\text{c}}k})(\lambda_{{\text{c}}j}-\lambda_{{\text{c}}k})p_{hi}^Rp_{hj}^R}, &i\neq j.
\end{cases}
\end{equation}
For simplicity, denoting $\p_r=(\p_R^{-1})^{\T}=(p_{ij}^r)_{r\times r}$, we have $\W_{\text{C}}^{-1}=\p_r\M_{\text{C}}^{-1}\p_r^{\T}$ and the specific elements are
\begin{equation*}
\begin{split}
\W_{\text{C}}^{-1}(i, j)&=\sum_{l=1}^r\sum_{k=1}^rp_{ik}^r\M_{\text{C}}^{-1}(k, l)p_{jl}^r\\
&=\sum_{l=1}^rp_{il}^r\M_{\text{C}}^{-1}(l, l)p_{jl}^r+\sum_{b=1,b\neq c}^r\sum_{c=1}^rp_{ic}^r\M_{\text{C}}^{-1}(c, b)p_{jb}^r\\
&=\sum_{l=1}^rp_{il}^rp_{jl}^r\frac{(1-\lambda_{{\text{c}}i}^2)\prod\limits_{k\neq i}(1-\lambda_{{\text{c}}i}\lambda_{{\text{c}}k})^2}{\prod\limits_{k\neq i}(\lambda_{{\text{c}}i}-\lambda_{{\text{c}}k})^2(p^{R}_{hi})^2}\\
&\quad +\sum_{b=1,b\neq v}^r\sum_{v=1}^rp_{iv}^rp_{jb}^r\frac{(1-\lambda_{{\text{c}}v}^2)(1-\lambda_{{\text{c}}b}^2)\prod\limits_{k< l}\prod\limits_{l=2}^{n}(1-\lambda_{{\text{c}}k}\lambda_{{\text{c}}l})}{(\lambda_{{\text{c}}v}-\lambda_{{\text{c}}b})^2\prod\limits_{k\neq v,b}(\lambda_{{\text{c}}v}-\lambda_{{\text{c}}k})(\lambda_{{\text{c}}b}-\lambda_{{\text{c}}k})p_{hv}^Rp_{hb}^R}.
\end{split}
\end{equation*}
%and{\tiny
%\begin{equation*}
%\begin{split}
%(\W_{\text{C}}^{-1})^2(i, j)&=\sum_{s=1}^r\W_{\text{C}}^{-1}(i, s)\W_{\text{C}}^{-1}(s, j)\\
%&=\sum_{s=1}^r\left(\sum_{l=1}^rp_{il}^rp_{sl}^r\frac{(1-\lambda_{{\text{c}}i}^2)\prod\limits_{k\neq i}(1-\lambda_{{\text{c}}i}\lambda_{{\text{c}}k})^2}{\prod\limits_{k\neq i}(\lambda_{{\text{c}}i}-\lambda_{{\text{c}}k})^2(p^{R}_{hi})^2}\right.\\
%&\quad\left. +\sum_{b=1,b\neq v}^r\sum_{v=1}^rp_{iv}^rp_{sb}^r\frac{(1-\lambda_{{\text{c}}v}^2)(1-\lambda_{{\text{c}}b}^2)\prod\limits_{k< l}\prod\limits_{l=2}^{n}(1-\lambda_{{\text{c}}k}\lambda_{{\text{c}}l})}{(\lambda_{{\text{c}}v}-\lambda_{{\text{c}}b})^2\prod\limits_{k\neq v,b}(\lambda_{{\text{c}}v}-\lambda_{{\text{c}}k})(\lambda_{{\text{c}}b}-\lambda_{{\text{c}}k})p_{hv}^Rp_{hb}^R}\right)\\
%&\left(\sum_{l=1}^rp_{sl}^rp_{jl}^r\frac{(1-\lambda_{{\text{c}}s}^2)\prod\limits_{k\neq s}(1-\lambda_{{\text{c}}s}\lambda_{{\text{c}}k})^2}{\prod\limits_{k\neq s}(\lambda_{{\text{c}}s}-\lambda_{{\text{c}}k})^2(p^{R}_{hs})^2}\right.\\
%&\quad\left. +\sum_{b=1,b\neq v}^r\sum_{v=1}^rp_{sv}^rp_{jb}^r\frac{(1-\lambda_{{\text{c}}v}^2)(1-\lambda_{{\text{c}}b}^2)\prod\limits_{k< l}\prod\limits_{l=2}^{n}(1-\lambda_{{\text{c}}k}\lambda_{{\text{c}}l})}{(\lambda_{{\text{c}}v}-\lambda_{{\text{c}}b})^2\prod\limits_{k\neq v,b}(\lambda_{{\text{c}}v}-\lambda_{{\text{c}}k})(\lambda_{{\text{c}}b}-\lambda_{{\text{c}}k})p_{hv}^Rp_{hb}^R}\right).
%\end{split}
%\end{equation*}}
Based on that, we can derive $\underline{\alpha}$ and $\underline{\beta}$ as
\begin{equation*}
\begin{split}
\underline{\alpha}=\sum_{i=1}^r\sum_{j=1}^r&\left(\sum_{l=1}^rp_{il}^rp_{jl}^r\frac{(1-\lambda_{{\text{c}}i}^2)\prod\limits_{k\neq i}(1-\lambda_{{\text{c}}i}\lambda_{{\text{c}}k})^2}{\prod\limits_{k\neq i}(\lambda_{{\text{c}}i}-\lambda_{{\text{c}}k})^2(p^{R}_{hi})^2}\right.\\
&\quad \left.+\sum_{b=1,b\neq v}^r\sum_{v=1}^rp_{iv}^rp_{jb}^r\frac{(1-\lambda_{{\text{c}}v}^2)(1-\lambda_{{\text{c}}b}^2)\prod\limits_{k< l}\prod\limits_{l=2}^{n}(1-\lambda_{{\text{c}}k}\lambda_{{\text{c}}l})}{(\lambda_{{\text{c}}v}-\lambda_{{\text{c}}b})^2\prod\limits_{k\neq v,b}(\lambda_{{\text{c}}v}-\lambda_{{\text{c}}k})(\lambda_{{\text{c}}b}-\lambda_{{\text{c}}k})p_{hv}^Rp_{hb}^R}\right)^2
\end{split}
\end{equation*}
\begin{equation*}
\begin{split}
\underline{\beta}&=\sum_{i=1}^r\sum_{j=1}^r
\left(
\sum_{s=1}^r\left(\sum_{l=1}^rp_{il}^rp_{sl}^r\frac{(1-\lambda_{{\text{c}}i}^2)\prod\limits_{k\neq i}(1-\lambda_{{\text{c}}i}\lambda_{{\text{c}}k})^2}{\prod\limits_{k\neq i}(\lambda_{{\text{c}}i}-\lambda_{{\text{c}}k})^2(p^{R}_{hi})^2}\right.\right.\\
&\quad\left. +\sum_{b=1,b\neq v}^r\sum_{v=1}^rp_{iv}^rp_{sb}^r\frac{(1-\lambda_{{\text{c}}v}^2)(1-\lambda_{{\text{c}}b}^2)\prod\limits_{k< l}\prod\limits_{l=2}^{n}(1-\lambda_{{\text{c}}k}\lambda_{{\text{c}}l})}{(\lambda_{{\text{c}}v}-\lambda_{{\text{c}}b})^2\prod\limits_{k\neq v,b}(\lambda_{{\text{c}}v}-\lambda_{{\text{c}}k})(\lambda_{{\text{c}}b}-\lambda_{{\text{c}}k})p_{hv}^Rp_{hb}^R}\right)\\
&\left(\sum_{l=1}^rp_{sl}^rp_{jl}^r\frac{(1-\lambda_{{\text{c}}s}^2)\prod\limits_{k\neq s}(1-\lambda_{{\text{c}}s}\lambda_{{\text{c}}k})^2}{\prod\limits_{k\neq s}(\lambda_{{\text{c}}s}-\lambda_{{\text{c}}k})^2(p^{R}_{hs})^2}\right.\\
&\quad\left.\left. +\sum_{b=1,b\neq v}^r\sum_{v=1}^rp_{sv}^rp_{jb}^r\frac{(1-\lambda_{{\text{c}}v}^2)(1-\lambda_{{\text{c}}b}^2)\prod\limits_{k< l}\prod\limits_{l=2}^{n}(1-\lambda_{{\text{c}}k}\lambda_{{\text{c}}l})}{(\lambda_{{\text{c}}v}-\lambda_{{\text{c}}b})^2\prod\limits_{k\neq v,b}(\lambda_{{\text{c}}v}-\lambda_{{\text{c}}k})(\lambda_{{\text{c}}b}-\lambda_{{\text{c}}k})p_{hv}^Rp_{hb}^R}\right)
\right)^2.
\end{split}
\end{equation*}

When $|\lambda_{{\text{c}}1}|<1$ and other eigenvalues satisfy $|\lambda_{{\text{c}}1}|\leq\cdots\leq|\lambda_{{\text{c}}\mu}|<1$, $|\lambda_{{\text{c}} \,\mu+1}|=\cdots=|\lambda_{{\text{c}} \,\mu+k}|=1$, $1<|\lambda_{{\text{c}} \,\mu+k+1}|\leq\cdots\leq|\lambda_{{\text{c}}r}|$, similar to the section~\ref{1drivecom}, the elements of the first $\mu$ rows and the first $\mu$ columns dominate.
Therefore, for $\W_{\text{C}}^{-1}(i, j)=\sum_{l=1}^r\sum_{k=1}^rp_{ik}^r\M_{\text{C}}^{-1}(k, l)p_{jl}^r$, we employ $\sum_{l=1}^\mu\sum_{k=1}^\mu p_{ik}^r\M_{\text{C}}^{-1}(k, l)p_{jl}^r$ to approximate it.
Accordingly, $\underline{\alpha}$ and $\underline{\beta}$ are
\begin{equation}\label{Aeq13}
\begin{split}
\underline{\alpha}=\sum_{i=1}^r\sum_{j=1}^r&\left(\sum_{l=1}^\mu p_{il}^rp_{jl}^r\frac{(1-\lambda_{{\text{c}}i}^2)\prod\limits_{k\neq i}(1-\lambda_{{\text{c}}i}\lambda_{{\text{c}}k})^2}{\prod\limits_{k\neq i}(\lambda_{{\text{c}}i}-\lambda_{{\text{c}}k})^2(p^{R}_{hi})^2}\right.\\
&\quad \left.+\sum_{b=1,b\neq v}^\mu\sum_{v=1}^\mu p_{iv}^rp_{jb}^r\frac{(1-\lambda_{{\text{c}}v}^2)(1-\lambda_{{\text{c}}b}^2)\prod\limits_{k< l}\prod\limits_{l=2}^{n}(1-\lambda_{{\text{c}}k}\lambda_{{\text{c}}l})}{(\lambda_{{\text{c}}v}-\lambda_{{\text{c}}b})^2\prod\limits_{k\neq v,b}(\lambda_{{\text{c}}v}-\lambda_{{\text{c}}k})(\lambda_{{\text{c}}b}-\lambda_{{\text{c}}k})p_{hv}^Rp_{hb}^R}\right)^2
\end{split}
\end{equation}
\begin{equation}\label{Aeq14}
\begin{split}
\underline{\beta}&=\sum_{i=1}^r\sum_{j=1}^r
\left(
\sum_{s=1}^r\left(\sum_{l=1}^\mu p_{il}^rp_{sl}^r\frac{(1-\lambda_{{\text{c}}i}^2)\prod\limits_{k\neq i}(1-\lambda_{{\text{c}}i}\lambda_{{\text{c}}k})^2}{\prod\limits_{k\neq i}(\lambda_{{\text{c}}i}-\lambda_{{\text{c}}k})^2(p^{R}_{hi})^2}\right.\right.\\
&\quad\left. +\sum_{b=1,b\neq v}^\mu\sum_{v=1}^\mu p_{iv}^rp_{sb}^r\frac{(1-\lambda_{{\text{c}}v}^2)(1-\lambda_{{\text{c}}b}^2)\prod\limits_{k< l}\prod\limits_{l=2}^{n}(1-\lambda_{{\text{c}}k}\lambda_{{\text{c}}l})}{(\lambda_{{\text{c}}v}-\lambda_{{\text{c}}b})^2\prod\limits_{k\neq v,b}(\lambda_{{\text{c}}v}-\lambda_{{\text{c}}k})(\lambda_{{\text{c}}b}-\lambda_{{\text{c}}k})p_{hv}^Rp_{hb}^R}\right)\\
&\left(\sum_{l=1}^\mu p_{sl}^rp_{jl}^r\frac{(1-\lambda_{{\text{c}}s}^2)\prod\limits_{k\neq s}(1-\lambda_{{\text{c}}s}\lambda_{{\text{c}}k})^2}{\prod\limits_{k\neq s}(\lambda_{{\text{c}}s}-\lambda_{{\text{c}}k})^2(p^{R}_{hs})^2}\right.\\
&\quad\left.\left. +\sum_{b=1,b\neq v}^\mu\sum_{v=1}^\mu p_{sv}^rp_{jb}^r\frac{(1-\lambda_{{\text{c}}v}^2)(1-\lambda_{{\text{c}}b}^2)\prod\limits_{k< l}\prod\limits_{l=2}^{n}(1-\lambda_{{\text{c}}k}\lambda_{{\text{c}}l})}{(\lambda_{{\text{c}}v}-\lambda_{{\text{c}}b})^2\prod\limits_{k\neq v,b}(\lambda_{{\text{c}}v}-\lambda_{{\text{c}}k})(\lambda_{{\text{c}}b}-\lambda_{{\text{c}}k})p_{hv}^Rp_{hb}^R}\right)
\right)^2.
\end{split}
\end{equation}
When $|\lambda_{{\text{c}}1}|=1$, we have $\underline{\alpha}\sim {\tau_f}^{-2}$ and $\underline{\alpha}\sim {\tau_f}^{-4}$, and
when $|\lambda_{{\text{c}}1}|>1$, we have
$\underline{\alpha}\thicksim \lambda_{{\text{c}}1}^{-4{\tau_f}}$, \text{and} $ \underline{\beta}\thicksim \lambda_{{\text{c}}1}^{-8{\tau_f}}.$

In the case of $m$ driver nodes with $m\leq r$, when $|\lambda_{{\text{c}}r}|<1$, compared with the case of $1$ driver node, we find only $q_{ij}^C, i, j=1, 2, \dots, r$ is different.
Denoting the indexes of driver nodes by $d_1, d_2, \cdots, d_m$,
then the corresponding matrix is $\B_r=[e_{d_1}, e_{d_2}, \cdots, e_{d_m}]\in \R^{r\times m}$.
Therefore, we have $\overline{\alpha}$ and $\overline{\beta}$ as
\begin{equation}\label{Aeq15}
\overline{\alpha}=\sum_{i=1}^r\sum_{j=1}^r\left( \sum_{k=1}^r\sum_{l=1}^r p_{ik}^Rp_{jl}^Rq_{kl}^C\frac{1}{1-\lambda_{{\text{c}}k}\lambda_{{\text{c}}l}}\right)^2
\end{equation}
and
\begin{equation}\label{Aeq16}
\begin{split}
\overline{\beta}&=\sum_{i=1}^r\sum_{j=1}^r\left( \sum_{s=1}^r\left(\sum_{k=1}^r\sum_{l=1}^r p_{ik}^Rp_{sl}^Rq_{kl}^C\frac{1}{1-\lambda_{{\text{c}}k}\lambda_{{\text{c}}l}}\right)
\left( \sum_{k=1}^r\sum_{l=1}^r p_{sk}^Rp_{jl}^Rq_{kl}^C\frac{1}{1-\lambda_{{\text{c}}k}\lambda_{{\text{c}}l}} \right) \right)^2
\end{split}
\end{equation}
with $q_{ij}^C=\sum_{k=1}^mp_{d_k\,i}^Rp_{d_k\,j}^R$.
Other cases are similar to the section~\ref{1drivecom} and thus omitted here.

\end{document}